\documentclass[10pt,a4paper,final]{iopart}

\usepackage{iopams}
\expandafter\let\csname equation*\endcsname\relax
\expandafter\let\csname endequation*\endcsname\relax
\usepackage{amsmath}
\usepackage{graphicx}
\usepackage[breaklinks=true,colorlinks=true,linkcolor=blue,urlcolor=blue,citecolor=blue]{hyperref}
\usepackage{fancyhdr, amssymb, amsthm}
\usepackage{accents, mathrsfs, tensor}
\usepackage{makecell, array}
\usepackage{pst-node}
\usepackage{texdraw}
\usepackage{tikz}
\usepackage{upgreek,xspace}
\usepackage{comment}
\usepackage{cite}
\usepackage{color, soul}
\sethlcolor{orange}

\usepackage{lipsum}
\makeatletter
\def\footnoterule{\relax%
  \kern15pt
  \hbox to \columnwidth{\vrule width 0.5\columnwidth height 0.4pt}
  \kern3pt}
\makeatother

\theoremstyle{definition}

\begin{document}

\def\Nu{\ensuremath\upnu\xspace}
\def\theequation{\arabic{section}.\arabic{equation}}
\newcommand{\pp}{\partial}
\newcommand{\aar}{\bar{a}}
\newcommand{\bb}{\beta}
\newcommand{\gm}{\gamma}
\newcommand{\Gm}{\Gamma}
\newcommand{\en}{\epsilon}
\newcommand{\ven}{\varepsilon}
\newcommand{\dd}{\delta}
\newcommand{\sg}{\sigma}
\newcommand{\kp}{\kappa}
\newcommand{\ld}{\lambda}
\newcommand{\oa}{\omega}
\newcommand{\hf}{\frac{1}{2}}
\newcommand{\be}{\begin{equation}}
\newcommand{\ee}{\end{equation}}
\newcommand{\bea}{\begin{eqnarray}}
\newcommand{\eea}{\end{eqnarray}}
\newcommand{\bse}{\begin{subequations}}
\newcommand{\ese}{\end{subequations}}
\newcommand{\nn}{\nonumber}
\newcommand{\bR}{\bar{R}}
\newcommand{\bP}{\bar{\Phi}}
\newcommand{\bS}{\bar{S}}
\newcommand{\bu}{{\boldsymbol u}}
\newcommand{\bt}{{\boldsymbol t}}
\newcommand{\bm}{{\boldsymbol m}}
\newcommand{\boa}{{\boldsymbol \omega}}
\newcommand{\bet}{{\boldsymbol \eta}}
\newcommand{\bW}{\bar{W}}
\newcommand{\vf}{\varphi}
\newcommand{\sn}{{\rm sn}}
\newcommand{\wh}{\widehat}
\newcommand{\ol}{\overline}
\newcommand{\wt}{\widetilde}
\newcommand{\ut}{\undertilde}
\newcommand*\xbar[1]{\hbox{\vbox{\hrule height 0.4pt\kern0.5ex\hbox{\kern-0.3em\ensuremath{#1}\kern-0.1em}}}}
\newcommand{\mathcolorbox}[2]{\colorbox{#1}{$\displaystyle #2$}}

\title[Shami Alsallami, Jitse Niesen, and Frank Nijhoff]{Closed-form modified Hamiltonians for integrable numerical integration schemes}

\author{S.A.M. Alsallami, J. Niesen, and F.W. Nijhoff}
\address{School of Mathematics, University of Leeds, Leeds LS2 9JT, United Kingdom}
\eads{mmsaa@leeds.ac.uk, j.niesen@leeds.ac.uk and f.w.nijhoff@leeds.ac.uk}

\date{\today}

\begin{abstract}

Modified Hamiltonians are used in the field of geometric numerical integration to show that symplectic schemes for Hamiltonian systems are accurate over long times. For nonlinear systems the series defining the modified Hamiltonian usually diverges. In contrast, this paper constructs and analyzes explicit examples of nonlinear systems where the modified Hamiltonian has a closed-form expression and hence converges. These systems arise from the theory of discrete integrable systems. We present cases of one- and two-degrees symplectic mappings arising as reductions of nonlinear integrable lattice equations, for which the modified Hamiltonians can be computed in closed form. These modified Hamiltonians are also given as power series in the time step by Yoshida's method based on the Baker-Campbell-Hausdorff series. Another example displays an implicit dependence on the time step which could be of relevance to certain implicit schemes in numerical analysis. In the light of these examples, the potential importance of integrable mappings to the field of geometric numerical integration is discussed.

\end{abstract}


\section{Introduction}

This paper combines ideas from discrete integrable systems and geometric numerical integration. Geometric numerical integrators are numerical methods that preserve geometric properties of the flow of a differential equation. They have relatively small error propagation in long time integrations, even for simple integration algorithms, cf.~\cite{hairer2006geometric}. Backward error analysis is the key to explaining this phenomenon. Backward error analysis states that the numerical scheme can be viewed as the time step of a modification of the original system. It is a well-known and essential fact that for a Hamiltonian system the modified system is Hamiltonian if the numerical scheme is symplectic. One then speaks of a modified Hamiltonian (MH) or an interpolating Hamiltonian. This MH is a Hamiltonian, the flow of which interpolates the iterations of the numerical scheme, i.e.\ its flow coincides with the computed points. Moser~\cite{moser1968lectures} developed a formal scheme for constructing the MH as an expansion and Benettin \& Giorgilli~\cite{benettin1994hamiltonian} established that this MH exists as an asymptotic series in the step size of the numerical scheme. This series is generally divergent, the necessary truncation of the series induces an error, which can be made exponentially small in the step size~\cite{hairer2006geometric, reich1999backward}. However, no error needs to be induced when the expansion for the MH is convergent. Examples of the MH for the harmonic oscillator and discussions of this situation can be found in~\cite{sanz1992symplectic, sanz1994numerical, skeel2001practical}. However, for nonlinear systems the convergence of the MH can only be achieved in exceptional circumstances, e.g.\ when the numerical scheme is integrable.

In the one-degree-of-freedom case, an autonomous Hamiltonian system always has a conserved quantity since the Hamiltonian itself is a first integral. For the discrete-time system this is no longer true and it is exceptional to have an invariant, which only happens if the system is integrable~\cite{mcmillan1971problem, quispel1988integrable, quispel1989integrable, veselov1991integrable}. A non-integrable map, however, is not expected to possess a globally defined invariant function on its phase space, and since the MH is an invariant of the numerical scheme by construction for such a map the MH can not exist as proper function (i.e.\ it can not have a convergent MH). On the other hand, if the numerical scheme is integrable, there must be a link between the invariant and the MH, possibly through a transcendental relation, cf.~\cite{field2003note}. In the multiple-degrees-of-freedom case, the situation is more subtle; for complete integrability the system needs to possess as many independent invariants as there are degree-of-freedom. This is the situation we explore in this paper with regard to the two-degrees-of-freedom case. However, to have a convergent MH it might not be necessary to have more than one first integral and the system may need to be only partially integrable (such a situation corresponds to so-called quasi-integrable systems). In this paper, we will not explore the latter possibility but consider the numerical scheme viewed as a dynamical map $(\boldsymbol{q},\boldsymbol{p})\longrightarrow(\overline{\boldsymbol{q}},\overline{\boldsymbol{p}})$ which is symplectic and completely integrable, i.e.\ possesses a full set of invariants~$\{I_j\}$ with~$I_j(\overline{\boldsymbol{q}}, \overline{\boldsymbol{p}})=I_j(\boldsymbol{q},\boldsymbol{p})$ which are independent and an involution with respect to the Poisson bracket, cf.~\cite{veselov1991integrable}. The integrable numerical schemes that we consider in this paper arise as reduction of nonlinear integrable lattice equations, which are integrable partial difference equations on a quadrilateral lattice, cf.~\cite{hietarinta2016discrete, nijhoff1995discrete, nijhoff1983direct}. These equations arise also as a numerical algorithms, e.g.~Pad{\'e} approximant and convergence acceleration algorithms~\cite{gragg1972thePadé}, and they are also important for the study of numerically induced chaos~\cite{taha1984}.

Reductions to maps are obtained from periodic initial value problems for such integrable lattice equations. We present examples of systems in one- and two-degrees-of-freedom arising from nonlinear integrable lattice equations, which give mappings with  invariants, and more than this, allow us to write the MHs {\it in closed form}. In particular, these examples arise from periodic initial value problems of the lattice versions of Korteweg-de Vries (KdV) and modified Korteweg-de Vries (MKdV) equations, they are constructed by updating the lattice variables along a diagonal shift so that the mapping is close to the identity mapping if the step size is small. This is convenient, because it allows us to obtain symplectic mappings from the consideration of an initial value problem on a two-dimensional lattice.

The outline of this paper is as follows. First, the MH due to Yoshida's method~\cite{yoshida1990construction} is reviewed in Section~\ref{sec. of Yoshida for MH}, then in Section~\ref{sec. of Existence MH}, the existence and expansion convergence of the MH are discussed. In Section~\ref{sec. of KdV and mKdV mappings}, we show how the lattice KdV and MKdV equations give rise to finite-dimensional mappings which do indeed carry a spectral interpretation, and that the invariants can be calculated systematically from the monodromy matrix constructed from the Lax pair. The Lax pair description of the mappings of KdV and MKdV types are then given. In Section~\ref{sec. of 1-degree Hamiltonian}, we study the one-degree-of-freedom Hamiltonian system arising from the lattice KdV and MKdV equations, and we show the derivation of the MH using Yoshida's method. We also write the MH in closed form expression using action-angle variables technique. The transition to multiple-degrees-of-freedom brings important new features, such as finite-gap integration technique~\cite{belokolos1994algebro}, which we explore in the rest of the paper, but to keep the discussion transparent, we restrict ourselves to two-degrees-of-freedom system. The multicomponent case requires the technique of separation of variables~\cite{sklyanin1995separation, van1976spectrum}, which we discuss in Section~\ref{sec. of SoV for KdV}, for the transformation to action-angle variables which involves the theory of genus-two abelian features, as developed in Section~\ref{sec. of 2-degrees Hamiltonian}. Finally, in Section~\ref{sec. of Conclusion}, we give a conclusion and discuss the scope for establishing a bridge between the theory of mappings and the mathematical structures in geometric integration.


\section{The modified Hamiltonian}
\setcounter{equation}{0}


\subsection{Yoshida's construction} \label{sec. of Yoshida for MH}

Symplectic integrators are numerical integration schemes for $N$-degrees-of-freedom Hamiltonian systems~$(\{q_j\}, \{p_j\}, H)$, which conserve the symplectic two-form $\sum_{j}{dq_j}\,{\wedge}\,{dp_j}$, where $j=1, \dots, N$.
In $2N$-dimensional phase space, the Poisson bracket is defined as
\be
\{F, G\}=\sum\limits_{j=1}^{N}\left(\frac{\partial F}{\partial q_j}\frac{\partial G}{\partial p_j}-\frac{\partial F}{\partial p_j}\frac{\partial G}{\partial q_j}\right),
\label{canonical Poisson bracket sec.2}
\ee
where $q_j,p_j$ are coordinates of the phase space, $j=1, \dots, N$.
Defining the differential operator $D_G$ by 
$$D_GF:=\{F, G\}\ ,$$ we have
$$\dot{z}_j=D_H\,z_j(q_j, p_j)=\{z_j, H\},\quad\textrm{where}\qquad z\equiv(q_1,\dots,q_N;p_1,\dots,p_N)\ .$$
Thus, the integrated $t$-flow of the equations of motion can be written by $$e^{tD_H}(z_j(q_j, p_j))\ .$$
Consider the symplectic Euler method
\be
	\overline{p}_j-p_j=-\tau\,\frac{\partial{H}(q_j, \overline{p}_j)}{\partial{q}_j},\qquad \overline{q}_j-q_j=\tau\,\frac{\partial{H}(q_j, \overline{p}_j)}{\partial{p_j}}\ ,
	\label{Euler method Yoshida construction}
\ee
applied to a Newtonian type $$H(q_1,\dots, q_N; p_1,\dots, p_N)=T(p_1,\dots, p_N)+V(q_1,\dots, q_N)\ .$$
Equations (\ref{Euler method Yoshida construction}) actually form a discrete analogue of the usual Hamilton's equations so we will refer to them as discrete Hamilton equations for short, and call the corresponding $H$ the discrete Hamiltonian. Thus, the discrete Hamiltonian is the generating function for the canonical transformation formation~(\ref{Euler method Yoshida construction}).
The symplectic Euler method is represented by
\be e^{\tau D_V}e^{\tau D_T}\ .\ee
The Baker-Campbell-Hausdorff (BCH) formula (cf.~\cite{oteo1991baker, varadarajan1974}) tells us that the product of two exponentials can be expressed as a single exponential, i.e. $$e^Xe^Y=e^Z\ ,$$ where $Z$ is given by the following infinite series of nested commutators:
$$Z=X+Y+\frac{1}{2}\,[X,Y]+\frac{1}{12}\,([X,[X,Y]]+[Y,[Y,X]])+\cdots\ .$$
Hence, the product $e^{\tau D_V}e^{\tau D_T}$
can be written as $e^{\tau D_{{H}^{\ast}}}$, where using $$[D_F, D_G]=D_{\{G,F\}}\ ,$$
we have
\begin{equation}
	{H}^{\ast}=T+V+\frac{\tau}{2}\,\{T,V\}+\frac{\tau^2}{12}\,\left(\{T,\{T,V\}\}+\{V,\{V,T\}\}\right)+\cdots\ .
	\label{eq:BCH for MH (general)}
\end{equation}
Thus, the symplectic Euler method follows the exact $\tau$ evolution of ${H}^{\ast}$. In other words, ${H}^{\ast}$ is the MH.
For the canonical Poisson bracket (\ref{canonical Poisson bracket sec.2}), this result can be written as
$${H}^{\ast}=H+\tau H_1+\tau^2H_2+\cdots$$
where
$$H_1=-\frac{1}{2}\sum_{j=1}^{N}\frac{\partial H}{\partial p_j}\frac{\partial H}{\partial q_j}, ~ H_2=\frac{1}{12}\sum_{i,j=1}^{N}\left(\frac{\partial H}{\partial p_j}\frac{\partial H}{\partial p_i}\frac{\partial^2 H}{\partial q_i\partial q_j}+\frac{\partial H}{\partial q_j}\frac{\partial H}{\partial q_i}\frac{\partial^2 H}{\partial p_i\partial p_j}\right), ~ \cdots\ .$$


\subsection{Existence and convergence} \label{sec. of Existence MH}

It is well known that for linear systems (quadratic Hamiltonians), the expansion for the MH is convergent. It used to be believed that the linear case is the only case of a discretization arising from a Hamiltonian of the form $$H=\hf\,p^2+V(q)\ ,$$ which gives a convergent series for the MH, while for nonlinear systems, the BCH expansion for the MH does not converge, cf.~\cite{hairer2006geometric, sanz1994numerical, skeel2001practical}. If the series does not converge, the MH may not exist as a proper function and as such provides only a formal invariant.  However, the paper~\cite{field2003note} gives examples of systems coming from the theory of discrete integrable systems, for which there exists a closed-form expression for the MH, indicating that the BCH expansion should converge. Generally, if the BCH expansion converges, the MH defines an invariant for the mapping. This provides a link to the theory of discrete  integrable systems.

In this paper, we build on the work in~\cite{field2003note} and give further examples of discrete integrable systems which, when viewed as an application of the symplectic Euler method, have a closed-form expression for the MH. The corresponding Hamiltonian systems are associated with the interpolating flow of these integrable mappings. In particular, we give examples of one- and two-degrees-of-freedom systems which arise from nonlinear integrable lattice equations. The construction of these examples will be exposed in the next section.


\section{Integrable lattice equations and dynamical mappings}\label{sec. of KdV and mKdV mappings}
\setcounter{equation}{0}


\subsection{The lattice KdV equation}\label{sec. of KdV mapping and invariants}

\begin{figure}[ht]
\centering
\begin{tikzpicture}[every node/.style={minimum size=0.25cm}]
\tikzstyle{solid node}=[circle,draw,inner sep=1,fill=black]
\draw[thick,-] (0,2) -- (0,3);
\draw[thick,-] (2,2) -- (0,2) node[above right=-0.1] {$\mathrm{u}_0$};
\draw[thick,-] (2,0) -- (2,2) node[above right=-0.3] {$\mathrm{u}_1$};
\draw[thick,-] (2,0) -- (4,0) node[above right=-0.3] {$\mathrm{u}_3$};
\draw[thick,-] (4,0) -- (4,-2) node[above right=-0.3] {$\mathrm{u}_4$};
\draw[thick,-] (4,-2) -- (6,-2) node[above right=-0.3] {};
\draw[thick,-] (6,-3) -- (6,-2) node[above right=-0.3] {$\mathrm{u}_5$};
\draw[thick,-] (0,0) -- (2,0) node[above right=-0.3] {$\mathrm{u}_2$};
\draw[thick,-] (0,2) -- (0,0) node[below left=-0.3] {$\overline{\mathrm{u}}_1$};
\draw[thick,-] (4,-2) -- (4,-4);
\draw[thick,-] (5,-4) -- (4,-4) node[below left=-0.3] {$\overline{\mathrm{u}}_5$};
\draw[thick,-] (2,0) -- (2,-2);
\draw[thick,-] (4,-2) -- (2,-2) node[below left=-0.3] {$\overline{\mathrm{u}}_3$};
\draw[thick,-] (-1,2) -- (0,2);
\draw[thick,dashed,-] (0,-2) -- (0,0);
\draw[thick,dashed,-] (2,-2) -- (0,-2) node[below left=-0.3] {$\overline{\mathrm{u}}_2$};
\draw[thick,dashed,-] (0,0) -- (-2,0) node[below left=-0.3] {$\overline{\mathrm{u}}_0$};
\draw[thick,dashed,-] (-2,0) -- (-2,1);
\draw[thick,dashed,-] (2,-4) -- (2,-2);
\draw[thick,dashed,-] (4,-4) -- (2,-4) node[below left=-0.3] {$\overline{\mathrm{u}}_4$};
\draw[thick,dashed,-] (4,-5) -- (4,-4);
\node[solid node] at (-2,0) {};
\node[solid node] at (0,0) {};
\node[solid node] at (2,0) {};
\node[solid node] at (4,0) {};
\node[solid node] at (6,0) {};
\node[solid node] at (0,2) {};
\node[solid node] at (2,2) {};
\node[solid node] at (4,2) {};
\node[solid node] at (-2,-2) {};
\node[solid node] at (0,-2) {};
\node[solid node] at (2,-2) {};
\node[solid node] at (4,-2) {};
\node[solid node] at (6,-2) {};
\node[solid node] at (0,-4) {};
\node[solid node] at (2,-4) {};
\node[solid node] at (4,-4) {};
\node[solid node] at (2,-6) {};
\draw[thick,-] (7,-4) -- (7,-5);
\draw[thick,-] (6,-5) -- (7,-5) node[above left=-0.3] {$\mathrm{u}_{2P}=\mathrm{u}_0$};
\draw[thick,-] (7,-5) -- (9,-5) node[above right=-0.3] {$\mathrm{u}_1=\wt{\mathrm{u}}_0$};
\draw[thick,-] (7,-5) -- (7,-7) node[below left=-0.3] {$\overline{\mathrm{u}}_1=\wh{\mathrm{u}}_0$};
\draw[thick,-] (9,-5) -- (9,-7);
\draw[thick,-] (7,-7) -- (9,-7) node[above right=-0.3] {$\mathrm{u}_2=\wh{\wt{\mathrm{u}}}_0$};
\draw[thick,-] (9,-7) -- (9,-8);
\draw[thick,-] (9,-7) -- (10,-7);
\node[solid node] at (5,-7) {};
\node[solid node] at (7,-9) {};
\node[solid node] at (5,-9) {};
\node[solid node] at (7,-5) {};
\node[solid node] at (7,-7) {};
\node[solid node] at (9,-3) {};
\node[solid node] at (9,-5) {};
\node[solid node] at (9,-7) {};
\end{tikzpicture}
	\caption{Standard staircase of periodic initial data on lattice}
\label{Lattice figure}
\end{figure}
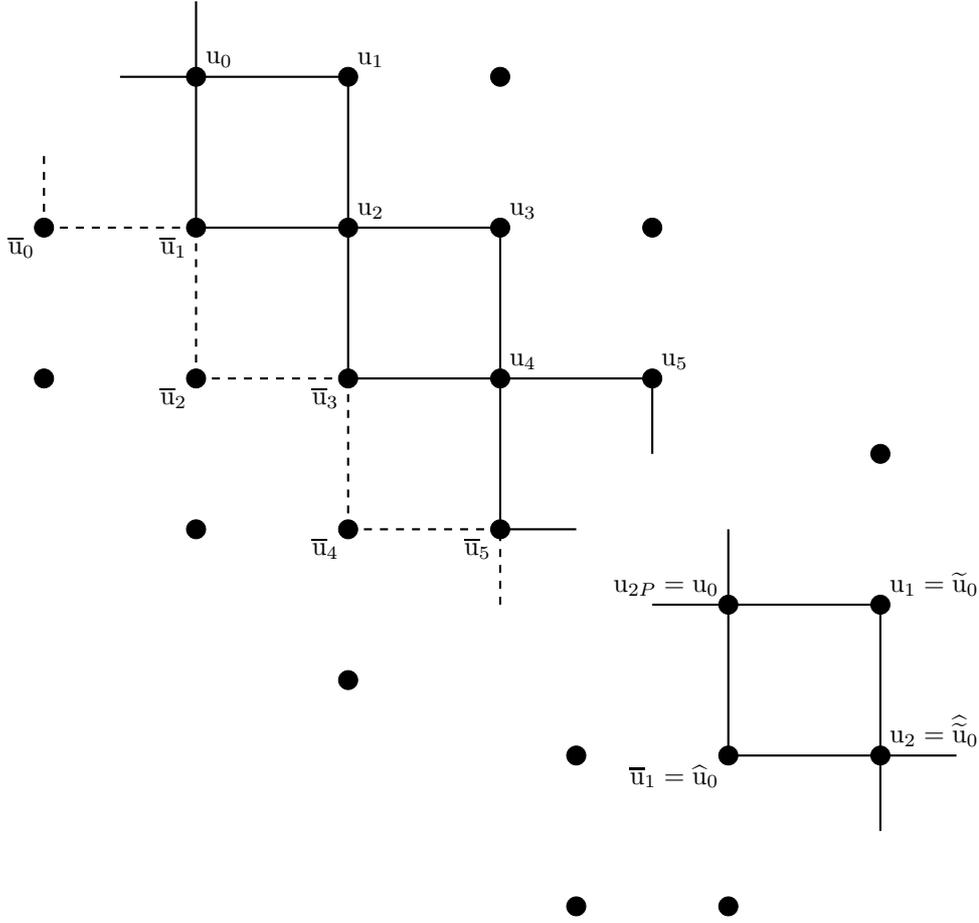

The lattice version of the KdV equation that we prefer to work with is the following nonlinear partial difference equation~\cite{hirota1977-1981, nijhoff1983direct, papageorgiou1990integrable},
\begin{equation}
(p-q+\widehat{u}-\widetilde{u})\, (p+q-\widehat{\widetilde{u}}+u)=p^2-q^2
\label{KdV eq}\ .
\end{equation}
Here $u:=u(n, m)$ is the dynamical variable at the lattice site $(n, m)$ with $n, m\in \mathbb Z$, and $\wt{\phantom{a}}$ and $\wh{\phantom{a}}$ are shorthand notations for translations on the lattice, i.e.~$\widetilde{u}:=u(n+1,m)$ and $\widehat{u}:=u(n,m+1)$, as shown in figure~\ref{Lattice figure}. Furthermore, $p$ and $q$ are complex-valued lattice parameters.
Equation (\ref{KdV eq}) arises as the compatibility condition of a pair of linear problems (Lax pair) defining the shifts of two component vector functions $\Phi(k)$ in the $n$ and $m$ directions,
\begin{equation}
	(p-k)\widetilde{\Phi}(k)=\mathcal{L}(k)\Phi(k), \ \quad   \ (q-k)\widehat{\Phi}(k)=\mathcal{M}(k)\Phi(k)\ ,
	\label{Lax pair KdV}
\end{equation}
where $\mathcal{L}(k)$ is given by
\begin{equation}
\mathcal{L}(k)=\begin{pmatrix}
		p-\widetilde{u} & 1\\
		k^2-p^2+(p-\widetilde{u})(p+u) & p+u
	\end{pmatrix},
\label{lax matrix for KdV}
\end{equation}
and where $\mathcal{M}(k)$ is given by a similar matrix obtained from (\ref{lax matrix for KdV}) by making the replacements $p \rightarrow q$ and $\wt{\phantom{a}}\rightarrow\wh{\phantom{a}}$. The parameter $k$ is the spectral parameter. An important feature of the equation (\ref{KdV eq}) is that it arises from a discrete action principle. The action for the KdV lattice equation~(\ref{KdV eq}) reads
\be
\mathcal{S}=\sum_{n,m\in \mathbb Z}\left\lbrack u_{n,m}\left({u}_{n+1,m}-{u}_{n,m+1}\right)+\epsilon\,\delta\,\textrm{log}\,(\epsilon+u_{n,m}-{u}_{n+1,m+1})\right\rbrack ,
\label{eq:action for KdV}
\ee
in which $\delta=p-q$, $\epsilon=p+q$. We note that we will be concerned with equation~(\ref{eq:action for KdV}) in section~\ref{sec. of 2-degrees Hamiltonian}.

Let us now consider initial value problems for (\ref{KdV eq}) on the lattice. One way of doing this is to give initial data on a horizontal line which leads to a nonlocal scheme, cf.~\cite{wiersma1987lattice}. In this paper, we are concerned with another type of initial value problem, which gives rise to a local iteration scheme: we assign initial data on a staircase on the lattice, cf.\ ref.~\cite{papageorgiou1990integrable}. By a staircase we mean a sequence of neighbouring lattice sites with $m$ and $n$ nondecreasing, as e.g.\ illustrated in figure~\ref{Lattice figure}. From the fact that equation~(\ref{KdV eq}) at each site involves only the four variables situated on the four lattice sites around a simple plaquette, it follows that the information on these staircases evolves diagonally through the lattice along parallel staircases. Furthermore, because of the convexity of the staircase configuration, the initial-value problem is well-posed. Although staircases of variable length and height stairsteps can be considered, for the sake of clarity we use a standard staircase of an even-periodic configuration of initial data.
Thus, we choose initial data on the standard staircase through the origin $(n, m)=(0,0)$, as depicted in figure~\ref{Lattice figure}, namely
	$$u(j,j)=:\mathrm{u}_{2j},\qquad u(j+1,j)=:\mathrm{u}_{{2j}+1}\qquad (j\in \mathbb Z)\ .$$
In fact, more general initial value configuration could be considered as well, but we will not do so here.
In this paper, we take a different point of view from the one expounded in ref.~\cite{papageorgiou1990integrable}. We perform iterations by updating the lattice variables $u$ along a diagonal shift rather than vertical shift; that is we define
	$$u(j-1,j+1)=:\overline{\mathrm{u}}_{2j},\qquad u(j,j+1)=:\overline{\mathrm{u}}_{{2j}+1},$$
	using the lattice KdV (\ref{KdV eq}). In this way we obtain the following mapping,
	\begin{equation}
		\overline{\mathrm{u}}_{2j}=\mathrm{u}_{2j}-\delta+\frac{\epsilon\,\delta}{\epsilon-\ol{\mathrm{u}}_{2j+1}+\ol{\mathrm{u}}_{2j-1}}, \quad
		\overline{\mathrm{u}}_{{2j}+1}=\mathrm{u}_{{2j}+1}-\delta+\frac{\epsilon\,\delta}{\epsilon-\mathrm{u}_{2j+2}+\mathrm{u}_{2j}}\ .
		\label{eq:KdV map}
	\end{equation}
We choose a diagonal shift so that the mapping is close to the identity mapping if $\delta$ is small.
By introducing now the differences
	$$X_j:= \mathrm{u}_{{2j}+1}-\mathrm{u}_{{2j}-1},\qquad Y_j:= \mathrm{u}_{{2j}+2}-\mathrm{u}_{2j} \qquad (j\in \mathbb Z)\ ,$$ equations (\ref{eq:KdV map}) can be reduced to the rational mapping
\begin{equation}
		\overline{X}_j=X_j+\frac{\epsilon\,\delta}{\epsilon-Y_j}-\frac{\epsilon\,\delta}{\epsilon-Y_{j-1}}, \quad
		\overline{Y}_j=Y_j-\frac{\epsilon\,\delta}{\epsilon-\overline{X}_j}+\frac{\epsilon\,\delta}{\epsilon-\overline{X}_{j+1}}
		\label{eq:reduced KdV map}\ .
	\end{equation}
The mapping (\ref{eq:reduced KdV map}) arises as the compatibility condition of a linear (Zakharov-Shabat type of) problem, that is easily obtained from the linear problem for the lattice equation (\ref{KdV eq}). By using a special property of the Lax matrices $\mathcal{L}$ and $\mathcal{M}$ of (\ref{Lax pair KdV}), (\ref{lax matrix for KdV}), it turns out that we can perform at each site of the staircase a gauge transformation such that we obtain expressions for the $L$ and $M$ matrices in terms of the reduced variables $X_j$ and $Y_j$ only. In this way one obtains the representation
\be
\ol{L}_j\,M_j=M_{j+1}\,L_j
\label{compatibility condition KdV map}
\ee
for the mapping (\ref{eq:reduced KdV map}), in which
\bse
\begin{align}
{L}_j&=\begin{pmatrix}
		y_{j} & 1\\
		\ld & 0
	\end{pmatrix}\begin{pmatrix}
		x_{j} & 1\\
		\ld-\epsilon\,\delta & 0
	\end{pmatrix} ,
\label{Lax matrix L for KdV map} \\
{M}_j&=\begin{pmatrix}
		{-\epsilon\,\delta}/\ol{x}_j & 1\\
		\ld & -\ol{x}_j
	\end{pmatrix}\begin{pmatrix}
		x_{j}-{\epsilon\,\delta}/{y_{j}} & 1\\
		\ld & 0
	\end{pmatrix} ,
	\label{Lax matrix M for KdV map}
\end{align}\ese
where we have used the abbreviations $x_{j}\equiv\epsilon-X_j$, $y_{j}\equiv\epsilon-Y_j$ and $\ld\equiv k^2-q^2$.
For convenience a detailed description is outlined in \ref{appendix:a}.

We impose now the even periodicity condition $\mathrm{u}_{2(j+P)}=\mathrm{u}_{2j},~ \mathrm{u}_{2(j+P)+1}=\mathrm{u}_{2j+1}$. Note that $P\,(P=2, 3, ...)$ can be interpreted as the period along the two diagonals of the lattice corresponding to the staircase. It is easy to see that these periodic conditions are compatible with the lattice equation, and hence will be preserved after iteration of the mapping. This means that we have to supply (\ref{eq:reduced KdV map}) with the periodicity constraints,
\begin{equation}
\sum_{j=1}^{P}X_j=0,\qquad \sum_{j=0}^{P-1}Y_j=0\ .
\label{periodicity constraints for KdV}
\end{equation}
Equation\ (\ref{eq:reduced KdV map}) for $j= 1, 2, ..., P$ together with the constraint (\ref{periodicity constraints for KdV}) is a $2(P-1)$-dimensional integrable mapping. The simplest case $P=2$ corresponds to the mapping,
\begin{equation}
		\overline{X}=X+\frac{2\,\epsilon\,\delta\, Y}{\epsilon^2-Y^2}, \qquad
		\overline{Y}=Y-\frac{2\,\epsilon\,\delta\,\overline{X}}{\epsilon^2-\overline{X}^2}
		\label{KdV map in X,Y (g=1)}\ .
	\end{equation}
For the case of $P=3$, the corresponding mapping reads in terms of four variables $X_1, X_2, Y_1, Y_2$,
\bse\label{mapping of KdV (P=3)in X and Y}
\begin{align}
		&\overline{X}_1=X_1+\frac{\epsilon\,\delta}{\epsilon-Y_1}-\frac{\epsilon\,\delta}{\epsilon+Y_1+Y_2}\ , \,\, \quad
		\overline{Y}_1=Y_1+\frac{\epsilon\,\delta}{\epsilon-\overline{X}_2}-\frac{\epsilon\,\delta}{\epsilon-\overline{X}_1}\qquad\, , \\
		&\overline{X}_2=X_2+\frac{\epsilon\,\delta}{\epsilon-Y_2}-\frac{\epsilon\,\delta}{\epsilon-Y_1} \qquad, \   \
		\  \ \overline{Y}_2=Y_2+\frac{\epsilon\,\delta}{\epsilon+\overline{X}_1+\overline{X}_2}-\frac{\epsilon\,\delta}{\epsilon-\overline{X}_2}\ .
\end{align}\ese

All mappings obtained in this way exhibit $P-1$ nontrivial integrals which can be found in a straightforward way by exploiting the Lax representation (\ref{Lax pair KdV}) of the original lattice KdV equation.
In order to do this, we need to define the monodromy matrix $\mathcal{T}(k)$ as
\begin{equation}
	\mathcal{T}(k)=\prod_{v=0}^{\substack{\curvearrowleft\\N-1}}\mathcal{L}(\mathrm{u}_{v+1},\mathrm{u}_v)
	\label{monodromy matrix for KdV sec.}\ ,
\end{equation}
where the $\curvearrowleft$ indicates that the factors in the product are arranged from right to left, and in which the translation matrices $\mathcal{L}(\mathrm{u}_{v+1},\mathrm{u}_v)$ represent either the Lax matrix in the $n$-direction ($v$ even), i.e.~$\mathcal{L}(k)$, or the Lax matrix in the $m$-direction ($v$ odd), i.e.~$\mathcal{M}(k)$. They are of the form
\begin{equation}
\mathcal{L}(\mathrm{u}_{v+1},\mathrm{u}_v)=\begin{pmatrix}
	p_{v+1}-\mathrm{u}_{v+1} & 1\\
	k^2-p_{v+1}(\mathrm{u}_{v+1}-\mathrm{u}_v)-\mathrm{u}_{v+1}\mathrm{u}_v & p_{v+1}+\mathrm{u}_{v}
\end{pmatrix} ,
\label{Lax matrix for n,m-direction}
\end{equation}
where $p_{v+1}=p$ if $v$ is even and $p_{v+1}=q$ if $v$ is odd.
This leads to, cf.~\cite{papageorgiou1990integrable},
\begin{equation}
\begin{split}
	\textrm{tr}\,\mathcal{T}(k)=&\,\prod_{v=1}^{N}(p_v+p_{v+1}+\mathrm{u}_{v-1}-\mathrm{u}_{v+1})\\
	&+\sum_{j=1}^{N}(k^2-p_j^2)\prod_{\substack{{v=1}\\{v\neq j-1,j}}}^{N}(p_v+p_{v+1}+\mathrm{u}_{v-1}-\mathrm{u}_{v+1})\\
	&+\sum_{i<j=1}^{N}(k^2-p_i^2)(k^2-p_j^2)\prod_{\substack{{v=1}\\{v\neq i-1,i,j-1,j}}}^{N}(p_v+p_{v+1}+\mathrm{u}_{v-1}-\mathrm{u}_{v+1})\\
	&+\cdots .
	\end{split}
	\label{trace of matrix (general)}
\end{equation}
The coefficients of the powers of $k$ are the integrals of mappings.

Alternatively, having obtained the linear system in Zakharov-Shabat form (\ref{compatibility condition KdV map}), one can also construct the monodromy matrix $T(\lambda)$ by gluing the elementary translation matrices $L_j(\lambda)$ along the staircase over one period $P$, leading to
\be
{T}(\lambda):=\prod_{j=0}^{\substack{\curvearrowleft\\ {P-1}}}{L}_j(\lambda)
\label{monodrmy matrix using Lax pair (KdV Sec3)}\ .
\ee
The monodromy matrix $T(\lambda)$ from~(\ref{monodrmy matrix using Lax pair (KdV Sec3)}) is the same as $\mathcal{T}(k)$ from~(\ref{monodromy matrix for KdV sec.}), apart from a similarity transformation corresponding to a gauge transformation at the beginning and end point of the chain from $0$ to $2P$ (which by periodicity is performed by the same multiplying matrix), and the trace of the monodromy matrix leads to the same result as before.

In the case of $N=2P=4$, equation\ (\ref{trace of matrix (general)}) yields the integral of the mapping
\begin{equation}
	\mathcal{I}=X^2\,Y^2-\epsilon^2\,X^2-\epsilon^2\,Y^2-2\,\epsilon\,\delta\,X\,Y
	\label{eq:invariant(X,Y):KdV g=1}\ ,
\end{equation}
which is an elliptic curve that can be parametrized in terms of Jacobi elliptic functions, thus leading to explicit solutions to the corresponding mapping.
For $N=2P=6$, the invariants are calculated from (\ref{trace of matrix (general)}) as
\bse \label{eq:invariants(X,Y):KdV g=2}
\begin{align}\nonumber
\mathcal{I}_1=&\,\frac{1}{8}\,{\delta}^{2}\,{\epsilon}^{2} \left\lbrack {x_1} \left(3\,{y_0}+3\,{y_1}-{y_2}\right)+{x_2} \left(3\,{y_1}+3\,{y_2}-{y_0}\right)+{x_3}\left(3\,{y_0}+3\,{y_2}-{y_1} \right)\right\rbrack\\ \nonumber
&+\frac{1}{2}\,\delta\,\epsilon\left\lbrack {x_1}\,{x_2}\,{y_1}\left({y_0}-{y_2}\right)+{x_1}\,{x_3}\,{y_0}\left({y_2}-{y_1} \right)+{x_2}\,{x_3}\,{y_2} \left({y_1}-{y_0} \right)\right\rbrack\\
&+{x_1}\,{x_2}\,{x_3}\,{y_0}\,{y_1}\,{y_2}\ , \\[0.1cm] \nonumber
\mathcal{I}_2=&\,\delta\,\epsilon \left\lbrack {x_1} \left({y_0}-{y_1} \right)+{x_2}\left( {y_1}-{y_2} \right)+{x_3}\left( {y_2}-{y_0} \right) \right\rbrack\\
&+{x_1}\,{x_2}\,{y_1}\left({y_0}+{y_2} \right)+{x_1}\,{x_3}\,{y_0}\left( {y_1}+{y_2} \right)+{x_2}\,{x_3}\,{y_2}\left( {y_0}+{y_1} \right) ,
\end{align}\ese
in which
$$x_j\equiv\epsilon-X_j,\quad y_j\equiv\epsilon-Y_j\quad (j=1,2),\quad y_0\equiv\epsilon+Y_1+Y_2, \quad x_3\equiv\epsilon+X_1+X_2\ .$$

\subsection{The lattice MKdV equation} \label{sec. of mKdV mapping and invariants}

Let us now consider the lattice version of the MKdV equation
\begin{equation}
pv\widehat{v}+q\widehat{v}\widehat{\widetilde{v}}=qv\widetilde{v}+p\widetilde{v}\widehat{\widetilde{v}}
\label{mKdV eq}\ .
\end{equation}
The parameters $p,q$ denote as before the lattice parameters, and the notations for the translations in the lattice direction are as before in the lattice KdV case.
Equation\ (\ref{mKdV eq}) arises as the compatibility condition of the linear system,
\begin{equation*}
	(p-k)\widetilde{\Psi}(k)=\mathfrak{L}(k)\Psi(k)\ ,\quad (q-k)\widehat{\Psi}(k)=\mathfrak{M}(k)\Psi(k)\ ,
	\label{Lax pair MKdV}
\end{equation*}
where $\mathfrak{L}(k)$ and $\mathfrak{M}(k)$ are given by
$$\mathfrak{L}(k)=
\left(\begin{array}{cc}
 p & \widetilde{v}\\
k^2/v & p\,\widetilde{v}/v \\
\end{array} \right), \  \ \mathfrak{M}(k)=\left(\begin{array}{cc}
q & \widehat{v}\\
k^2/v & q\,\widehat{v}/v \\
\end{array} \right) .
$$
Equation (\ref{mKdV eq}) is related to the lattice KdV equation via a Miura transformation as follows
\begin{equation}
	p-q+\widehat{u}-\widetilde{u}=\frac{p\,\widetilde{v}-q\,\widehat{v}}{v}
	\label{Miura transformation from KdV to mKdV}\ .
\end{equation}
On the level of the linear system, this reflects a gauge transformation of the form
$$k^2\,\Psi(k)=\mathcal{U}\,\Phi(k),\quad \mathcal{U}=\begin{pmatrix}
		s & v\\
		k^2 & 0
	\end{pmatrix},$$
in which $s=(p-\widetilde{u})\,v-p\,\widetilde{v}$ and $\widetilde{s}=p\,v-(p+u)\,\widetilde{v}$, leading to $$ \label{L(mkdv) given by L(kdv)} \mathfrak{L}(k)=\widetilde{\mathcal{U}}\,\mathcal{L}\,\mathcal{U}^{-1},\quad \widetilde{\mathcal{U}}\,F\,\mathcal{U}^{-1}=\frac{\widetilde{v}}{k^2}\,E\ ,\quad
\textrm{where}\quad E=\begin{pmatrix}
		0 & 1\\
		0 & 0
	\end{pmatrix},~ ~F=\begin{pmatrix}
		0 & 0\\
		1 & 0
	\end{pmatrix} .$$
It should be noted that a similar relation with the same $s$ holds for $\mathfrak{M}(k)$, by just replacing $p \rightarrow q$ and $\wt{\phantom{a}}\rightarrow\wh{\phantom{a}}$\ , i.e.\ $$\mathfrak{M}(k)=\widehat{\mathcal{U}}\,\mathcal{M}\,{\mathcal{U}}^{-1}.$$

Let us now consider initial value problems for (\ref{mKdV eq}) on the lattice in precisely the same way as before with initial data on staircases as the one depicted in figure~\ref{Lattice figure}, namely $$v(j,j)=:\mathrm{v}_{2j},\qquad v(j+1,j)=:\mathrm{v}_{{2j}+1}\qquad(j\in \mathbb Z)\ .$$
We perform iterations by updating the lattice variables $v$ along a diagonal shift, i.e.
	$$v(j-1,j+1)=:\overline{\mathrm{v}}_{2j},\qquad v(j,j+1)=:\overline{\mathrm{v}}_{{2j}+1}\ ,$$
using the lattice MKdV (\ref{mKdV eq}). One obtains the mapping
	\begin{equation}
		\overline{\mathrm{v}}_{2j}=\mathrm{v}_{2j}\,\frac{\mathrm{v}_{2j-1}+\rho\,\mathrm{v}_{2j+1}}{\mathrm{v}_{2j+1}+\rho\,\mathrm{v}_{2j-1}} , \quad
		\overline{\mathrm{v}}_{2j+1}=\mathrm{v}_{2j+1}\,\frac{\mathrm{v}_{2j}+\rho\,\mathrm{v}_{2j+2}}{\mathrm{v}_{2j+2}+\rho\,\mathrm{v}_{2j}}
		\label{mKdV mapping with b'}\ ,
	\end{equation}
where $\rho=p/q$. Again, we use a diagonal shift so that the mapping is close to the identity mapping if $\rho$ is close to $1$. We can reduce the system~(\ref{mKdV mapping with b'}) in terms of the logarithmic variables $$X_j:= \textrm{log}\,\frac{\mathrm{v}_{2j+1}}{\mathrm{v}_{2j-1}},\qquad Y_j:= \textrm{log}\,\frac{\mathrm{v}_{2j+2}}{\mathrm{v}_{2j}} \qquad (j \in \mathbb Z)\ .$$ 
The reduced mapping turns out to be
\begin{equation}
		\overline{X}_{j}=X_{j}+\textrm{log}\, \frac{(\rho\,+e^{Y_{j-1}})(1+\rho\,e^{Y_{j}})}{(\rho\,+e^{Y_{j}})(1+\rho\,e^{Y_{j-1}})},\quad
		\overline{Y}_{j}=Y_{j}+\textrm{log}\, \frac{(\rho\,+e^{\overline{X}_j})(1+\rho\,e^{\overline{X}_{j+1}})}{(\rho\,+e^{\overline{X}_{j+1}})(1+\rho\,e^{\overline{X}_{j}})}
		\label{eq:reduced MKdV map}\ .
	\end{equation}
In a similar way as before we can find from the Zakharov-Shabat system for the lattice equation~(\ref{mKdV eq}), a linear system for the mapping in terms of reduced variables,
\bse\label{Lax matrices for MKdV maps}
\begin{align}
{L}_j&=\left(\begin{array}{cc} q & y_j\\ k^2 & q\,y_j\end{array} \right)
\left(\begin{array}{cc} p & x_j\\ k^2 & p\,x_j\end{array} \right),\label{Lax matrix L for MKdV maps} \\
{M}_j&=\begin{pmatrix} -p & {(1+\rho\,x_j)}/{(\rho+x_j)}\\ k^2\,\ol{x}_j^{-1} & -p\,\ol{x}_j^{-1}{(1+\rho\,x_j)}/{(\rho+x_j)}\end{pmatrix}
\begin{pmatrix} q & x_j{(1+\rho\,y_j)}/{(\rho+y_j)}\\ k^2 & q\,x_j{(1+\rho\,y_j)}/{(\rho+y_j)}\end{pmatrix} ,\label{Lax matrix M for MKdV maps}
\end{align}\ese
in which we have used the abbreviations $x_j\equiv e^{X_j}, y_j\equiv e^{Y_j}$.
The Zakharov-Shabat equations, i.e.\ $$\ol{L}_j\,M_j=M_{j+1}\,L_j\ ,$$
leads to the mapping (\ref{eq:reduced MKdV map}) as compatibility condition.
We again impose the periodicity condition
\begin{equation}
\sum_{j=1}^{P}X_j=0,\qquad \sum_{j=0}^{P-1}Y_j=0\ .
\label{periodicity constraints mKdV}
\end{equation}
The mapping (\ref{eq:reduced MKdV map}) together with the conditions (\ref{periodicity constraints mKdV}) is a $2(P-1)$-dimensional mapping that exhibits $P-1$ invariants, that can be constructed in the same way as in the case of the lattice KdV by using the monodromy matrix. The invariants can also be constructed from the invariants of the KdV by using the Miura transformation~(\ref{Miura transformation from KdV to mKdV}).

In the case of $N=4$, the corresponding mapping is given as
 \begin{equation}
 \begin{split}
		\overline{X}=X+2\,\textrm{log}\,\frac{1+\rho\,e^{Y}}{\rho+e^{Y}} , \qquad
		\overline{Y}=Y+2\,\textrm{log}\,\frac{\rho+e^{\overline{X}}}{1+\rho\,e^{\overline{X}}}\ .
		 \end{split}
		\label{MKdV map in X,Y (g=1)}
	\end{equation}
We find the invariant
\begin{equation}
	\mathcal{I}=e^{X-Y}+e^{Y-X}+2\,\rho\,(e^{X}+e^{Y}+e^{-X}+e^{-Y})+\rho^2\,(e^{X+Y}+e^{-(X+Y)})
	\label{eq:invariant(X,Y):MKdV g=1}\ ,
\end{equation}
which is an elliptic curve.
 
In both the KdV and MKdV case we have similar structures, i.e. mappings and invariants. In fact, both cases are examples of discrete dynamical systems that are covered by the discrete version of the Arnold-Liouville theorem as formulated by Veselov, cf.~\cite{veselov1991integrable}. It was established in the paper~\cite{capel1991complete} that the invariants of the mappings are in involution.
It is a hallmark of integrability that these invariants themselves generate commuting flows using the Poisson bracket structures which are compatible with the discrete mappings. The task in the remainder of the paper is to make a connection between the invariants and the relevant MHs.


\section{The modified Hamiltonian of one-degree-of-freedom}\label{sec. of 1-degree Hamiltonian}
\setcounter{equation}{0}

In this section, we focus on the one-degree-of-freedom case obtained from the simplest reduction, leaving the systems with two-degrees-of-freedom to the rest of the paper.
We present a first example arising from the KdV reduction in section \ref{sec. of KdV example(g=1)}. The MH of the symplectic mapping, as given by Yoshida's approach, is then given. The MH is also written in closed form by using action-angle variables. We present a second example coming from the MKdV case in section \ref{sec. of mKdV example(g=1)}. We also do the same for the MH as in the first example.


\subsection{The KdV map example}\label{sec. of KdV example(g=1)}

Consider the integrable mapping
\begin{equation}
		\overline{p}=p+\frac{2\, \epsilon\, \delta\, q}{\epsilon^2-q^2}, \qquad
		\overline{q}=q-\frac{2\, \epsilon\, \delta\, \overline{p}}{\epsilon^2-\overline{p}^2}
		\label{KdV map in q,p(g=1)}\ ,
	\end{equation}
which is the mapping (\ref{KdV map in X,Y (g=1)}) in which we identify $X:= p$ and $Y:= q$.
We have the standard Poisson brackets
	\begin{equation*}
\{q,q\}=\{p,p\}=0  ,\quad  \{q,p\}=1\ ,
\end{equation*}
which is preserved by the map (\ref{KdV map in q,p(g=1)}).
Indeed, the mapping (\ref{KdV map in q,p(g=1)}) is a canonical transformation with the generating function
\begin{equation}
	{H}(q, \overline{p})=\epsilon\, \delta\, \textrm{log}\,(\epsilon^2-\overline{p}^2)+\epsilon\, \delta\, \textrm{log}\,(\epsilon^2-q^2) 
	\label{Hamiltonian eq. for KdV case P=2)}\ ,
\end{equation}
through the equations
\begin{equation}
	\overline{p}-p=-\frac{\partial H}{\partial q}, \qquad \overline{q}-q=\frac{\partial H}{\partial \overline{p}}\ ,
	\label{discrete Hamilton's equations KdV g=1}
\end{equation}
which, once again, we can regard as discrete analogue of the Hamilton equations, where (\ref{Hamiltonian eq. for KdV case P=2)}) can be regarded as the discrete Hamiltonian. The usual approach in geometric integration is to start with a Hamiltonian differential equation and then create a map using the same Hamiltonian as in the differential equation. Even if the differential equation is integrable, the resulting map will generically be non-integrable. In contrast, here we start from an integrable map and interpret the generating function of the map as a discrete Hamiltonian. In that case, there is an underlying continuous Hamiltonian flow whose Hamiltonian is given by an invariant of the map which is different from the discrete Hamiltonian.
More specifically, the mapping (\ref{KdV map in q,p(g=1)}) conserves the quantity
\begin{equation}
	\mathcal{I}=p^2\,q^2-\epsilon^2\,(p^2+q^2)-2\, \epsilon\, \delta\,p\,q
	\label{eq:invariant(q,p):KdV g=1}\ ,
\end{equation}
and, therefore, by definition is an integrable map (i.e. it is a symplectic mapping with an invariant).
By considering $\delta$ in the mapping (\ref{KdV map in q,p(g=1)}) to be the step size (i.e.\ $\delta$ plays the rule of $\tau$ in the expansion~(\ref{eq:BCH for MH (general)})), we obtain the following expansion for the MH:
\begin{equation}
	\begin{split}
	{H}^{\ast}=&\,\delta\,\epsilon\,\textrm{log}\,(\epsilon^2-p^2)+\delta\,\epsilon\,\textrm{log}\,(\epsilon^2-q^2)-\frac{2\,\delta^2\,\epsilon^2\,p\,q}{(\epsilon^2-p^2)(\epsilon^2-q^2)}\\
	&-\frac{2\,\delta^3\,\epsilon^3\,(\epsilon^2\, p^2+\epsilon^2\,q^2+2\,p^2\,q^2)}{3\,(\epsilon^2-p^2)^2(\epsilon^2-q^2)^2}
	-\frac{4\,\delta^4\,\epsilon^4\,p\,q\,(\epsilon^2+p^2)(\epsilon^2+q^2)}{3\,(\epsilon^2-p^2)^3(\epsilon^2-q^2)^3}
	+O(\delta^5)\ .
	\end{split}
	\label{MH by BCH (KdV g=1)}
\end{equation}
At this stage, it is not obvious that this expansion converges, however, we will show that we actually have a closed-form expression for the MH. This will provide a connection between the MH and the invariant, taking into account that if the MH exists, then it must be expressed in terms of the invariant.
As the mapping is integrable, we can employ a transformation to action-angle variables to derive an interpolating Hamiltonian where the canonical momenta are the invariants of the map.

A change to action-angle variables is a canonical (symplectic) transformation to a new set of phase space coordinates, such that the new momenta are invariants of the system and the coordinates evolve in a linear fashion. This can be viewed as an application of the Hamilton-Jacobi method, cf.~\cite{goldstein1956classical}. A discussion of action-angle variables for integrable mappings can be found in refs.~\cite{arnold1989mathematical, bruschi1991integrable}.

In the one-degree-of-freedom system under consideration, the relevant canonical transformation $(q,p)\longrightarrow(\mathcal{Q},\mathcal{P})$ is given by means of a generating function, $S(q,\mathcal{P})$ as
\be
	p=\frac{\partial{S}\,(q,\mathcal{P})}{\partial{q}}, \qquad \mathcal{Q}=\frac{\partial{S}\,(q,\mathcal{P})}{\partial\mathcal{P}}\ ,
	\label{eq:action-angle defining eqs(q,p)KdV,g=1)}
\ee
with
\begin{equation}\label{eq:K (KdV, g=1)}
	K(\mathcal{Q},\mathcal{P})=H+\frac{\partial S}{\partial t}\ ,
\end{equation}
being the transformed Hamiltonian, noting that in the case of a transformation to action-angle variables the $K$ only depends on $\mathcal{P}$. Integrating the first system of equations, we get $S$ up to an arbitrary function of the invariant
\begin{equation}
	S(q,\mathcal{P})=\int_{q^0}^{q} \frac{\epsilon\,\delta\,q'+\sqrt{\epsilon^2\,\delta^2\,q'^2-(\epsilon^2-q'^2)(\epsilon^2\, q'^2+\mathcal{P})}}{q'^2-\epsilon^2}\, dq'\ ,
\label{eq:S function (KdV g=1)}
\end{equation}
and consequently we obtain
\begin{equation}
	\mathcal{Q}(q,\mathcal{P})=\int_{q^0}^{q}\frac{1}{2\,\sqrt{\delta^2\, \epsilon^2\,q'^2-(\epsilon^2-q'^2)(\epsilon^2\, q'^2+\mathcal{P})}}\,dq'\ .
	\label{eq:coordinate Q (KdV g=1)}
\end{equation}

In the context of the mapping (\ref{KdV map in q,p(g=1)}), the relevant continuous Hamiltonian flow is the one whose Hamiltonian is given by the invariant of the map for which we have the Hamilton's equations
\begin{equation}
	\dot p=-\frac{\partial \mathcal{I}}{\partial q}, \qquad \dot q=\frac{\partial \mathcal{I}}{\partial p}\ .
	\label{Hamilton eqs with H=I (KdV g=1)}
\end{equation}
This system actually defines an interpolating flow where the trajectory of the system (\ref{Hamilton eqs with H=I (KdV g=1)}) and the orbit of the map share the level set of the invariant. For the system (\ref{discrete Hamilton's equations KdV g=1}) the new momentum is defined to be $\mathcal{P}=\mathcal{I}$, hence we integrate (\ref{Hamilton eqs with H=I (KdV g=1)}) by quadrature to obtain
\begin{equation}
t=\int_{0}^{\mathcal{E}(t|\epsilon, \delta,\mathcal{P})}\frac{1}{2\,\sqrt{\delta^2\, \epsilon^2\,q^2-(\epsilon^2-q^2)(\epsilon^2\, q^2+\mathcal{P})}}\ dq\ ,
\label{t elliptic integral}
\end{equation}
which defines the relevant elliptic function $\mathcal{E}(t|\epsilon, \delta,\mathcal{P})$ in terms of an elliptic integral of the first kind~\cite{hancock2004lectures}.
In fact, the modified Hamiltonian ${H}^{\ast}$ coincides with the canonical transformed Hamiltonian $K$ obtained by applying the canonical transformation (\ref{eq:action-angle defining eqs(q,p)KdV,g=1)}), viewed as a function of $\mathcal{Q}, \mathcal{P}$. Hamilton's equations in the new variables imply
\be\label{eqs:Hamilton's eqs(P,Q)KdV, g=1}
\dot{\mathcal{P}}=-\frac{\partial {H}^{\ast}}{\partial \mathcal{Q}}=0,\qquad \dot{\mathcal{Q}}=\frac{\partial {H}^{\ast}}{\partial \mathcal{P}}=\Nu\ ,
\ee
which tell us on the one hand that ${H}^{\ast}$ is a function of $\mathcal{P}$ alone, and on the other hand that ${H}^{\ast}$ is obtained by integrating $\Nu$ with respect to $\mathcal{P}$ up to an arbitrary function of the invariant.

In order to apply the canonical transformation to the map (\ref{KdV map in q,p(g=1)}), we introduce the ``frequency" $\Nu$ as the discrete time-one step
\begin{equation}
	\Nu=\int_{q}^{\overline{q}}\frac{\partial{p}}{\partial{\mathcal{P}}}\, dq', \quad\textrm{so that}\quad \ol{\mathcal{Q}}-\mathcal{Q}=\Nu\ ,
	\label{eq:frequency-general(KdV g=1)}
\end{equation}
which crucially depends on $\mathcal{P}$ only.
Since the time flow $t$ interpolates the map (\ref{KdV map in q,p(g=1)}), the iteration of the map is a time-one step stroboscope of the continuous time, so the integral (\ref{t elliptic integral}) can be subdivided into uniform time-one iterate.
Thus, we can choose an initial point at $q$, and use the system (\ref{KdV map in q,p(g=1)}) to compute $\overline{q}$, starting from $q=0$ we then obtain
$$\overline{q}=\frac{2\,\epsilon^2\,\delta\,\sqrt{-\mathcal{P}}}{\epsilon^4+\mathcal{P}}\ .$$
Thus, the frequency (\ref{eq:frequency-general(KdV g=1)}) is given as
\begin{equation}
\Nu=\int_0^{\frac{2\,\epsilon^2\,\delta\,\sqrt{-\mathcal{P}}}{\epsilon^4+\mathcal{P}}}\frac{1}{2\,\sqrt{\delta^2\, \epsilon^2\,q^2-(\epsilon^2-q^2)(\epsilon^2\, q^2+\mathcal{P})}}\ dq\ ,
\label{eq:elliptic integral(KdV, g=1)}
\end{equation}
and hence
\begin{equation}
	{H}^{\ast}(\mathcal{P})=\int^{\mathcal{P}}\int_0^{\frac{2\,\epsilon^2\,\delta\,\sqrt{-\mathcal{P}'}}{\epsilon^4+\mathcal{P}'}}\frac{1}{2\,\sqrt{\delta^2\, \epsilon^2\,q^2-(\epsilon^2-q^2)(\epsilon^2\, q^2+\mathcal{P}')}}\, dq\,d\mathcal{P}'
	\label{eq: closed form-(KdV g=1)}\ .
\end{equation}
The inside integral is a definite integral and the outside integral is an indefinite integral which is determined up to an integration constant. Equation (\ref{eq: closed form-(KdV g=1)}) is a closed-form expression for the MH of the map~(\ref{KdV map in q,p(g=1)}).

Writing equation\ (\ref{eq:elliptic integral(KdV, g=1)}) as a series in $\delta$ and fixing the integration constant such that for $\delta=0$, the integrating over $\mathcal{P}$ vanishes, we obtain
\begin{equation}
{H}^{\ast}(\mathcal{P})=\delta\, \epsilon\, \textrm{log}\,(\mathcal{P}+\epsilon^4)+\frac{2\,\delta^3\,\epsilon^3\,\mathcal{P}}{3\,(\mathcal{P}+\epsilon^4)^2}-\frac{4\,\delta^5\,\epsilon^5\,\mathcal{P} \left(\mathcal{P}-2\,\epsilon^4\right)}{15\,(\mathcal{P}+\epsilon^4)^4} +O(\delta^7)\ .
\label{MH in ivariants( KdV g=1)}
	\end{equation}
Unsurprisingly, expansions (\ref{MH by BCH (KdV g=1)}) and (\ref{MH in ivariants( KdV g=1)}) are the same, where the integration constant in (\ref{MH in ivariants( KdV g=1)}) is fixed. This matching can be seen on two steps. Firstly, one needs to insert the invariant~$\mathcal{P}$ and expand~(\ref{MH in ivariants( KdV g=1)}) in orders of $\delta$. Secondly, one needs to rearrange the series arising from step one in orders of $\delta$, which can be done by combining all terms with the same order of $\delta$.


\subsection{The MKdV map example}\label{sec. of mKdV example(g=1)}

Consider the integrable mapping
 \begin{equation}
\overline{p}=p+2\,\textrm{log}\,\frac{1+\rho\,e^{q}}{\rho+e^{q}} , \qquad
		\overline{q}=q+2\,\textrm{log}\,\frac{\rho+e^{\overline{p}}}{1+\rho\,e^{\overline{p}}}
		\label{MKdV map in p,q (g=1)}\ ,
\end{equation}
which is the mapping (\ref{MKdV map in X,Y (g=1)}) where we identify $X:= p$ and $Y:= q$. We have the standard invariant symplectic structure
\begin{equation*}
d\overline{q} \wedge d\overline{p} =dq \wedge dp\ ,
\end{equation*}
which is preserved by the map~(\ref{MKdV map in p,q (g=1)}).
The mapping~(\ref{MKdV map in p,q (g=1)}) is in fact a canonical transformation, and the generating function $H$ (Hamiltonian) of the mapping is in this case found to be
\begin{equation}
	H(q, \overline{p})=2\int_0^{\overline{p}}\textrm{log}\,\frac{\rho+e^{\xi}}{1+\rho\,e^{\xi}}\,d\xi+2\int_0^q\textrm{log}\,\frac{\rho+e^{\xi}}{1+\rho\,e^{\xi}}\,d\xi
	\label{Hamiltonian eq. for mKdV case P=2}\ .
\end{equation}
Note that this Hamiltonian can be written in terms of dilogarithm functions using the well-known integral representation.
The discrete-time Hamilton equations are written as
\begin{equation*}
	\overline{p}-p=-\frac{\partial{H}}{\partial q},\qquad \overline{q}-q=\frac{\partial H}{\partial \overline{p}}
	\label{Hamilton equations mKdV}\ ,
\end{equation*}
and the mapping (\ref{MKdV map in p,q (g=1)}) conserves the function
\begin{equation}
	\mathcal{I}=e^{p-q}+e^{q-p}+2\,\rho\,(e^{p}+e^{q}+e^{-p}+e^{-q})+\rho^2\,(e^{p+q}+e^{-(p+q)})
	\label{eq:invariant(q,p):MKdV g=1}\ .
\end{equation}

In this example we actually have a different type of situation from what we had in the previous example since the step size in the mapping~(\ref{MKdV map in p,q (g=1)}) is much more hidden (implicit).
Setting $\tau=\rho-1$ and using~(\ref{eq:BCH for MH (general)}) the expansion for the MH when applied to~(\ref{Hamiltonian eq. for mKdV case P=2}) can be written as a series in orders of $\tau$ as the following:
\be\begin{split}
{H}^{\ast}=&\,2\,\tau\left[ p+q-2\,\log\,(1+{{\rm e}^{p}})(1+{{\rm e}^{q}})\right]\\
	&-{\tau}^{2}\left[p+q-2\,\log\,(1+{{\rm e}^{p}})(1+{{\rm e}^{q}})+ {\frac{4\,({{\rm e}^{p+q}}+1)}{(1+{{\rm e}^{p}})(1+{{\rm e}^{q}})}}\right] \\
    &+\frac{2\,{\tau}^{3}}{3}\left[{\frac{3\,({{\rm e}^{p+2q}}+{{\rm e}^{2p+q}}+{{\rm e}^{p}}+{{\rm e}^{q}})}{(1+{{\rm e}^{p}})^2(1+{{\rm e}^{q}})^2}}+\frac{2\,(3\,{{\rm e}^{2(p+q)}}+8\,{{\rm e}^{p+q}}+3)} {(1+{{\rm e}^{p}})^2(1+{{\rm e}^{q}})^2}\right]\\
    &+\frac{2\,{\tau}^{3}}{3}\left[p+q-2\,\log\,(1+{{\rm e}^{p}})(1+{{\rm e}^{q}}) \right]
    +O(\tau^4)\ .
\label{MH by BCH (MKdV g=1)}
\end{split}
\ee
Once again, it is not obvious that this expansion converges, however, finding the connection between the MH and the invariant will essentially assert that there is a closed-form expression for the MH.
Again, since the mapping (\ref{MKdV map in p,q (g=1)}) is integrable we can follow the action-angle prescription of section~\ref{sec. of KdV example(g=1)} to derive the MH. Indeed, we can define the invariant~(\ref{eq:invariant(q,p):MKdV g=1}) as the new momentum $\mathcal{P}$ and its canonically conjugated variable as the new coordinate $\mathcal{Q}$. The canonical transformation from $({q},{p})$ to $(\mathcal{Q},\mathcal{P})$ can be parametrized in terms of a generating function $S(q,\mathcal{P})$ as
\begin{equation*}
	p=\frac{\partial{S}(q,\mathcal{P})}{\partial{q}}, \qquad \mathcal{Q}=\frac{\partial{S}(q,\mathcal{P})}{\partial{\mathcal{P}}}\ .
\end{equation*}
Similarly as in example \ref{sec. of KdV example(g=1)}, the corresponding frequency $\Nu$ is given by
\begin{equation}
\Nu=\int_{q}^{\ol{q}}{\frac {dq'}{\sqrt{\left(2\,\rho^{2}+\mathcal{P}+2\right)\left(\mathcal{P}-2\,\rho^2-4\,\rho\,e^{q'}-4\,\rho\,e^{-q'}-2\right)}}}\ .
\label{eq:frequency mKdV}
\end{equation}
Choosing an initial point at $q$ starting from $q=0$, one then obtains
$$\ol{q}=2\,\textrm{log}\left(\frac{\rho+\Delta}{1+\rho\,\Delta}\right),$$
in which we define the shorthand
$$
\Delta\equiv \frac{\mathcal{P}-4\,\rho-\sqrt {\left(2\,\rho^2+\mathcal{P}+2\right)\left(\mathcal{P}-2\,{\rho}^{2}-8\,\rho-2\right)}}{2\,(1+\rho)^2}.
$$
The integral~(\ref{eq:frequency mKdV}) is also an elliptic integral of the first kind.
By using the Hamilton's equations with the action-angle coordinates, and denoting the new Hamiltonian by ${H}^{\ast}$, we obtain
\begin{equation}
{H}^{\ast}(\mathcal{P})=\int^{\mathcal{P}} \int_{0}^{2\,\textrm{log}\left(\frac{\rho+\Delta}{1+\rho\,\Delta}\right)}{\frac {dq\, d\mathcal{P}'}{\sqrt{\left(2\rho^{2}+\mathcal{P}'+2\right)\left(\mathcal{P}'-2\rho^2-4\rho\,e^{q}-4\rho\,e^{-q}-2\right)}}}\ ,
\label{closed form for mKdV P=2}
\end{equation}
up to an integration constant.
Equation~(\ref{closed form for mKdV P=2}) is a closed-form expression for the MH of the map~(\ref{MKdV map in p,q (g=1)}).
Again, setting $\tau=\rho-1$ and fixing the constant of integration over $\mathcal{P}$ such that for $\tau=0$, the integral vanishes, we can write (\ref{closed form for mKdV P=2}) as a series in $\tau$ as follows
\begin{equation}
\begin{split}
	{H}^{\ast}(\mathcal{P})=&\,{-2}\,\tau\,{\textrm{log}\,{(\mathcal{P}+4)}}+\tau^2\left( \textrm{log}\,{(\mathcal{P}+4)}-\frac{8}{(\mathcal{P}+4)}\right)\\
	&-\frac{2\,\tau^3}{3} \left( \textrm{log}\,{(\mathcal{P}+4)}-\frac{4}{(\mathcal{P}+4)}-\frac{24}{(\mathcal{P}+4)^2}\right)
	+O(\tau^4)\ .
\end{split}
	\label{MH in invariants(MKdV g=1)}
\end{equation}
Equation~(\ref{MH in invariants(MKdV g=1)}) can once again be expanded in orders of $\tau$ and, after following the same procedure of section~\ref{sec. of KdV example(g=1)}, the matching with series~(\ref{MH by BCH (MKdV g=1)}) will be seen.

The present examples of one-degree-of-freedom systems already display some useful techniques from Hamiltonian mechanics. As we shall observe next, the transition to higher-degrees-of-freedom systems requires some further technology from the theory of discrete integrable systems. We will restrict ourselves to two-degrees-of-freedom system to keep the discussion transparent. The higher-degrees-of-freedom case has to be treated by the method of separation of variables in order to arrive at the action-angle variables, and we restrict ourselves to the case of KdV mappings.


\section{Separation of variables and finite-gap integration}\label{sec. of SoV for KdV}
\setcounter{equation}{0}

In the two-degrees-of-freedom case the situation becomes more complicated, which requires some mathematical techniques, such as separation of variables~\cite{sklyanin1995separation, van1976spectrum} and finite-gap integration~\cite{belokolos1994algebro}\footnote[10]{We note that in this section we are going to apply these techniques for higher-degrees-of-freedom of genus $g$.}.
Instead of transforming $(\boldsymbol{q}, \boldsymbol{p})$ directly to the action-angle variables $(\boldsymbol{\mathcal{Q}},\boldsymbol{\mathcal{P}})$, here we look for a canonical transformation to separated variables $(\boldsymbol{q}, \boldsymbol{p})  \longrightarrow   (\boldsymbol{\mu}, \boldsymbol{\eta})$ and a further transformation $(\boldsymbol{\mu},\boldsymbol{\eta})\longrightarrow (\boldsymbol{\mathcal{Q}},\boldsymbol{\mathcal{P}})$, combined with a description of the dynamics of the separated variables $(\boldsymbol{\mu},\boldsymbol{\eta})\longrightarrow(\ol{\boldsymbol{\mu}},\ol{\boldsymbol{\eta}})$.
We employ the key technique of separation of variables to the mappings of KdV type in section~\ref{sec. of SoV}. Then, we impose the Poisson bracket for the monodromy matrix, followed by the Poisson brackets between its entries, to establish the separation of variables transform as a canonical transformation, in section~\ref{sec. of PB for monodromy matrix and r-matrix}. We derive the discrete dynamics in terms of the separated variables (discrete Dubrovin equations) for this case in section~\ref{sec. of Dubrovin eqs}. We establish the continuous-time evolution for the auxiliary spectrum generated by the invariants in section~\ref{sec. of intrp. flow}. We introduce a canonical transformation to action-angle variables on the basis of structure in terms of the separation variables in section~\ref{sec. of Generating functions}.


\subsection{Separation of variables}\label{sec. of SoV}

The method of separation of variables (SoV) plays an important role in studying Liouville integrable systems. The SoV originated from the development of Hamiltonian mechanics as a method to separate an $N$-degrees-of-freedom system into one-degree-of-freedom system through the Hamilton-Jacobi equation for particular Hamiltonians. The SoV approach has been applied to many families of finite-dimensional integrable systems, cf.~refs.~\cite{arnold1989mathematical, sklyanin1995separation}.

Taking the monodromy matrix ${T}(\ld)$ in the form
\be
{T}(\ld)=\left( \begin{array}{cc} 
A(\ld) & B(\ld) \\ 
C(\ld) & D(\ld) 
\end{array} \right),
\label{monodrmy matrix with A,B,C,D (SoV Sec)}
\ee
it is well known that in the periodic problems for equations of KdV type, the roots $\{\mu_j\}$, $j=1,\dots, g$, of the polynomial $B(\ld)$ define the so-called {\it auxiliary spectrum}, cf.~\cite{toda2012theory}, and they play the role of separation variables, cf.~\cite{van1976spectrum}. These roots correspond to the poles of the Baker-Akhiezer\footnote[11]{This terminology appears in the theory of discrete integrable systems, in particular, in finite-gap integration.} (BA) function by which we mean an eigenfunction of the Lax representation normalized in such way that its analytic behaviour as a function of the spectral parameter $\ld$, is characterized through the singularity structure.

The monodromy matrix ${T}(\ld)$ for genus $g$ has a natural grading in terms of the spectral parameter $\ld$:
\begin{equation} \label{monodrmy matrix lambda expansion KdV} 
{T}(\ld)= \begin{pmatrix}
\ld^{g+1}+\ld^{g}\,A_{g}+\cdots + A_0 & \ld^{g}\,B_g+\ld^{g-1}\,B_{g-1}+ 
\cdots +B_0\\ 
\ld\left(\ld^{g}\,C_{g}+ \cdots +C_0\right) & 
\ld\left(\ld^{g}+\ld^{g-1}\,D_{g-1}+ \cdots +D_0\right)  
\end{pmatrix} .
\end{equation}
We can write the trace of the monodromy matrix as
\begin{equation}\label{trace for KdV type (SoV)}
{\rm tr}\,{T}(\ld)=\mathcal{I}_0+\sum_{j=1}^{g}\mathcal{I}_j\,\ld^{j}+2\,\ld^{g+1}\ , \quad \mathcal{I}_j=A_j+D_{j-1}\  .
\end{equation}
As noted in section~\ref{sec. of KdV and mKdV mappings}, the coefficients $\mathcal{I}_j$, $j=0,\dots,g-1$, are the invariants of the map, while the top coefficient $\mathcal{I}_g=A_g+D_{g-1}$ is a Casimir with respect to the natural Poisson algebra associated with the dynamical map.
The discriminant of the curve takes on the form:
\be\label{eq:discriminant(KdV)}
\begin{split}
R(\ld)=\left({\rm tr}\,{T}\right)^2-4\,\det (T)
&=\left(2\,\ld^{g+1}+\sum_{j=0}^{g} \mathcal{I}_j\,\ld^{j}\right)^2-4\,\ld^{g+1}\left(\ld-\omega\right)^{g+1} ,
\end{split}\ee
where $\omega=\epsilon\,\delta$.
The auxiliary spectrum $B(\ld)$ has the following factorisation
\begin{equation}\label{eq:auxiliary spectrum B}
B(\ld) = B_g\prod_{j=1}^g (\ld-\mu_j)\    , 
\end{equation} 
where $B_g$ is a Casimir.

The linear problem for the BA function is the eigenvector
\be
{T}(\ld)\,\phi(\ld)=\upeta(\ld)\,\phi(\ld)\ ,
\label{eq. of eigenvector}
\ee
of the monodromy matrix ${T}(\ld)$ corresponding to the eigenvalue $\upeta(\ld)$ of the spectral curve
\be
\det\,({T}(\ld)-\upeta)=0
\label{eq:spectral curve}\ ,
\ee
which defines a hyperelliptic curve of genus $g$.
This provides that a normalisation of the eigenvectors $\phi$ is fixed
\be
\overset{\rightarrow}{\alpha}\, \phi=1\ ,
\label{normalisation vector}
\ee
where $\overset{\rightarrow}{\alpha}$ is a row vector suitably chosen.
The pair $({\ld}, \upeta)$ can be thought of as a point of the spectral curve~(\ref{eq:spectral curve}).
The BA function $\phi$ is then a meromorphic function on the spectral curve.

It is easy to see that the pairs $({\mu}_j, {\eta}_j)$, where $\eta_j=\upeta(\mu_j)$, thus defined satisfy the separation equations~(\ref{eq:spectral curve}), which express the fact that the $({\mu}_j,{\eta}_j)$ are lying on the spectral curve.
The canonicity of the variables $({\mu}_j, {\eta}_j)$ should be verified independently.
No general recipe is known how to guess the proper (that
is producing canonical variables) normalisation for the BA function. In our case the following normalisation works,
\be
\overset{\rightarrow}{\alpha}=(1,0)\ .
\label{eq:normalisation for (KdV)}
\ee
From the linear equation (\ref{eq. of eigenvector}) and normalisation (\ref{normalisation vector}) we derive that $$\overset{\rightarrow}{\alpha}\,{T}\,\phi=\upeta\ ,$$ and hence,
\begin{equation}
\phi=\left( \begin{array}{c} 
\overset{\rightarrow}{\alpha} \\ 
\overset{\rightarrow}{\alpha}\, T(\ld)
\end{array} \right)^{-1}\left(\begin{array}{c} 
1 \\ 
\upeta  
\end{array} \right)
\label{eq:phi vector}.
\end{equation}
The following determinant has to vanish on the separation variables $\mu_j$,
\be
B(\ld)=\det \left( \begin{array}{c} 
\overset{\rightarrow}{\alpha} \\ 
\overset{\rightarrow}{\alpha}\, {T}(\ld) 
\end{array} \right)=0\ .
\label{eq:B(ld) with det}
\ee
The formula (\ref{eq:B(ld) with det}) for the separation variables appeared already in~\cite{kuznetsov1997separation}
in the case of standard normalisation $$\overset{\rightarrow}{\alpha}=\overset{\rightarrow}{\alpha}_0\equiv (0, 1)\ .$$
The rational functions $A(\ld)$ and $D(\ld)$ of the entries of the monodromy matrix ${T}({\ld})$ satisfy the following relations
\be
A(\mu_j)=\upeta^{+}(\mu_j)=\eta_j\quad \textrm{and}\quad D(\mu_j)=\upeta^{-}(\mu_j)=\eta'_j \ .
\label{eq: A,D(mu) relations}\ee

What remains is to establish the SoV transform as a canonical transformation, i.e.\ to verify that the Poisson brackets between the separation variables, have the canonical structure.
To do this we use the information about the Poisson brackets between the entries of monodromy matrix ${T}({\ld})$ provided by the classical $r,s$-matrix structure.


\subsection{Poisson brackets and r-matrix structures}\label{sec. of PB for monodromy matrix and r-matrix}

Let us remind the reader of the Lax matrices of the mapping of KdV type considered in section~\ref{sec. of KdV and mKdV mappings}. These matrices depends on a discrete variable $n$ labelling the sites along a chain of length $P$, and are given by
\be
L_n= \left(\begin{array}{cc} \ld+x_{n}\,y_{n}-\omega & y_{n}\\ \ld\,x_{n} & \ld\end{array}\right).
\label{Lax matrix for r-matrix KdV}
\ee
The key object here is the monodromy matrix $T(\ld)$ obtained by gluing the elementary translation matrices $L_n$ along a line connecting the sites $1$ and $P+1$ over one period $P$, namely
\begin{equation}
{T}(\ld):=\prod_{n=1}^{\substack{\curvearrowleft\\ P}}{L}_n(\ld)
\label{eq:monodrmy matrix by Lax (SoV sec)}\ ,
	\end{equation}
which is essentially the monodromy matrix (\ref{monodrmy matrix using Lax pair (KdV Sec3)}).
The Poisson bracket for the monodromy matrix ${T}({\ld})$ in terms of the $r,s$-matrix structure follows from the discrete version of the non-ultralocal Poisson bracket structure~\cite{nijhoff1992integrable} and reads as
	\begin{equation}
	\{\overset{1}{T},\overset{2}{T}\}=r_{12}^{+}\,\overset{1}{T}\,\overset{2}{T}-\overset{2}{T}\,\overset{1}{T}\,r_{12}^{-}-\overset{1}{T}\, s_{12}^{+}\, \overset{2}{T}+\overset{2}{T}\, s_{12}^{-}\, \overset{1}{T}
	\label{Poisson structure for monodromy matrix}\ .
	\end{equation}
In~(\ref{Poisson structure for monodromy matrix}) the superscripts $1,2$ for the operator matrix ${T}$ denote the corresponding factor on which this $T$ acts (acting trivially on the other factors), i.e.
$$\overset{1}{T}:=T\otimes \mathbb{I} \quad\textrm{and}\quad \overset{2}{T}:=\mathbb{I}\otimes T\ .$$
We note that the proof of equation~(\ref{Poisson structure for monodromy matrix}) can be found in \ref{appendix:c}.

In the classical case the traces of powers of the monodromy matrix are invariant under the mapping as a consequence of the discrete-time evolution
\be
\overline{T}(\ld)=M(\ld)\,T(\ld)\,{M(\ld)}^{-1}, \qquad
{M}=\begin{pmatrix}
		\ld+\left({\omega}/{\ol{x}}\right)\left({\omega}/{y}-x\right) & -{\omega}/{\ol{x}}\\
		\ld\left(x-\ol{x}-{\omega}/{y}\right) & \ld
\end{pmatrix} ,
\label{map T (SoV)}
\ee
where $M$ is $M_1$, and the periodicity condition $M_{P+1}=M_1$. Thus, this leads to a sufficient number of invariants which are obtained by expanding the traces in powers of the spectral parameter $\ld$.
The dynamical map in terms of the monodromy matrix is preserved by the Poisson bracket as a consequence of the compatibility condition of a discrete-time ZS (Zakharov-Shabat) system (\ref{compatibility condition KdV map}). The involution property of the classical invariants, which was proven in~\cite{capel1991complete}, follows also from the Poisson bracket
\be\label{eq: involution property}
\{\textrm{tr}\,{T}({\ld}),\, \textrm{tr}\,{T}({\ld}')\}=0\ ,\ee
which in turn follows from (\ref{Poisson structure for monodromy matrix}).

The classical $r,s$-matrix structure for the mapping of KdV type is given by
\be
r_{12}^{-}=\frac{\mathscr{P}_{12}}{\ld_1-\ld_2}, \quad s_{12}^{-}=\frac{1}{\ld_1}\,E_1F_2, \quad s_{12}^{+}=\frac{1}{\ld_2}\,F_1E_2, \quad r_{12}^{+}=r_{12}^{-}+s_{12}^{+}-s_{12}^{-}\ ,
\label{r,s-matrix structures for KdV}
\ee
in which $\ld_\alpha=k_\alpha^2-q^2, \alpha=1,2$, and the permutation matrix $\mathscr{P}_{12}$ and the matrices $E$, $F$ are given by
$$\mathscr{P}_{12}=\begin{pmatrix} 1 & 0 & 0 & 0 \\ 0 & 0 & 1 & 0 \\ 0 & 1 & 0 & 0 \\ 0 & 0 & 0 & 1 \end{pmatrix}, \quad E=\begin{pmatrix} 0 & 1 \\ 0 & 0 \end{pmatrix}, \quad F=\begin{pmatrix} 0 & 0 \\ 1 & 0 \end{pmatrix}.$$
Using equations~(\ref{Poisson structure for monodromy matrix}) and (\ref{r,s-matrix structures for KdV}) we can extract the following Poisson brackets between the entries of the monodromy matrix~(\ref{monodrmy matrix with A,B,C,D (SoV Sec)}):
\bse\label{Poisson brackets (SoV)} \begin{align}
\{ A(\ld_1), A(\ld_2)\}&=\frac{1}{\ld_2}B(\ld_1)\,C(\ld_2)-\frac{1}{\ld_1}B(\ld_2)\,C(\ld_1)\ , \\
\{ A(\ld_1), B(\ld_2)\}&=\frac{A(\ld_2)\,B(\ld_1)-B(\ld_2)\,A(\ld_1)}{\ld_1-\ld_2}+\frac{1}{\ld_2}\,B(\ld_1)\,D(\ld_2)\ ,\label{bracket of A1B2} \\
\{ A(\ld_1), D(\ld_2)\}&=\frac{\ld_1\,B(\ld_1)\,C(\ld_2)-\ld_2\,B(\ld_2)\,C(\ld_1)}{\ld_1\left(\ld_1-\ld_2\right)}\ , \\
\{B(\ld_1), B(\ld_2)\}&= 0\ ,  \\
\{B(\ld_1), D(\ld_2)\}&=\frac{\ld_1\,B(\ld_1)\,D(\ld_2)-\ld_2\,D(\ld_1)\,B(\ld_2)}{\ld_1\left(\ld_1-\ld_2\right)}\ ,\label{bracket of B1D2} \\
\{D(\ld_1), D(\ld_2)\}&= 0\ .
\end{align} \ese
Equations (\ref{Poisson brackets (SoV)}) can be used to establish the canonicity 
of the separation variables:
\begin{equation}\label{eq:Poisson structure for SoV} 
\{ \mu_{i},\mu_{j}\}=\{\eta_{i},\eta_{j}\}=0, \qquad \{\mu_{i}, \eta_{j}\}=\dd_{ij}\eta_{j}\ ,
\end{equation}
using $\upeta(\mu_{j})=A(\mu_j)$.


\subsection{Discrete Dubrovin equations}\label{sec. of Dubrovin eqs}

The Dubrovin equations arise in the theory of finite-gap integration as the equations governing the dynamics of the auxiliary spectrum or equivalently of the poles of the Baker-Akhiezer function. In \cite{nijhoff1999integrable} the finite-gap integration of mapping reductions of the lattice KdV equation was considered, cf.~also~\cite{cao2012finite}, for complementary results.
As a byproduct difference analogues of the Dubrovin equations, cf.~ref.~\cite{nijhoff2000discrete}, were derived which form the equations of the discrete motion of the auxiliary spectrum under the KdV mappings.

Let us first present the equations of discrete motion for the diagonal evolution in terms of the auxiliary spectrum, which are different from the equations for the vertical evolution given in~\cite{nijhoff2000discrete}.
The discrete dynamics in terms of the separated variables follows from~(\ref{trace for KdV type (SoV)}), (\ref{eq:discriminant(KdV)}) and~(\ref{map T (SoV)}), and is given by a coupled system of set of first-order difference equations for the $\mu_j$, details of derivation are given in~\ref{appendix:d}, namely
\bse\label{eq:Dubrovin eq}
\begin{align} \nonumber
&\left\lbrack{\mathscr M}^{-1}\left(\boldsymbol{\kp}\sqrt{R({\boldsymbol{\mu}})}-\mathcal{I}_0\,{\boldsymbol e}\right)\right\rbrack_j
 +\left\lbrack\xbar{\mathscr M}^{-1}\left(\ol{\boldsymbol{\kp}}\sqrt{R(\ol{\boldsymbol{\mu}})}-\mathcal{I}_0\,{\boldsymbol e}\right)\right\rbrack_j \\
&=\frac{2\,\mathcal{I}_0\,S_{g-j}(\boldsymbol{\mu})}{(-1)^{j}\prod_{i=1}^g\mu_i}+2\,B_g\,\frac{\ol{x}}{\omega}\,(-1)^{g-j+1}\left[ S_{g-j+1}(\ol{\boldsymbol{\mu}})-S_{g-j+1}(\boldsymbol{\mu})\right],\label{eq:Dub eq1} \\ \nonumber
&\left[{\mathscr M}^{-1}\left(\boldsymbol{\kp}\sqrt{R({\boldsymbol{\mu}})}-\mathcal{I}_0\,{\boldsymbol e}\right)\right]_g -\left[\xbar{\mathscr M}^{-1}\left(\ol{\boldsymbol{\kp}}\sqrt{R(\ol{\boldsymbol{\mu}})}-\mathcal{I}_0\,{\boldsymbol e}\right)\right]_g \\
&=\frac{2\,\mathcal{I}_0}{(-1)^{g}\prod_{i=1}^g\mu_i}-2\,B_g\,\ol{x}+2\,C_g\,\frac{\omega}{\ol{x}}\ ,
\label{eq:Dub eq2}\end{align}\ese
where $j=1,\dots,g$.
Thus, the actual discrete Dubrovin equations comprise two expressions (\ref{eq:Dub eq1}), (\ref{eq:Dub eq2}), and coupled through the $\ol{x}$ which can be eliminated by combination of both.
In equation~(\ref{eq:Dubrovin eq}), $\ol{\phantom{a}}$ is the shift in the discrete dynamical variables, $\boldsymbol{\mu}=\left(\mu_1,\dots,\mu_g\right)^t$ denotes the vector with entries $\mu_j$, $\mathcal{I}_0$ denotes the invariant given in equation~(\ref{trace for KdV type (SoV)}), $\boldsymbol{e}=(1,1,\dots,1)^t$, $\mathscr{M}$ denotes the Vandermonde matrix
\be \label{eq: VanderMonde matrix}
\mathscr{M}= \left(\begin{array}{ccc}
\mu_1 & \cdots & \mu_1^g \\ 
\vdots & & \vdots \\ 
\mu_g & \cdots &\mu_g^g  
\end{array} \right) , \ee
and ${\boldsymbol S}(\boldsymbol{\mu})$ is the vector of symmetric products $S_k$ of its arguments, i.e.\
\[ S_k(\mu_1,\dots,\mu_g)\equiv \sum_{i_1<i_2<\dots<i_k}^{g} 
\mu_{i_1}\mu_{i_2}\dots \mu_{i_k},\qquad S_0(\mu_1,\dots,\mu_g)=1\ . \]
We note that the $\boldsymbol{\kp}$ in~(\ref{eq:Dubrovin eq}) denotes the sign ${\kp}=\pm$ corresponding to the choice of sheet of the Riemann surface, subject to the condition $\ol{{\kp}}={\kp}$. The latter is the case $\ol{\kp}=\undertilde{\widehat{\kp}}$ since the bar shift $\ol{\phantom{a}}$ is the composition of two shifts
$\kp \rightarrow \undertilde{\kp}=-\kp$ and $\kp \rightarrow \widehat{\kp}=-\kp$, each of which provoke a change of sheet of the Riemann surface~\cite{nijhoff2000discrete, nijhoff1999integrable}.

In the case of $g=1$ the discrete Dubrovin equations (\ref{eq:Dubrovin eq}) reduce to set of two coupled equations, namely
\bse\label{Dubrovin eq.(coupled) g=1}
\begin{align}\label{eq1:Dubrovin eq. g=1}
&\frac{1}{\mu}\sqrt{R(\mu)}+\frac{1}{\ol{\mu}}\sqrt{R(\ol{\mu})}
=\mathcal{I}_0\left(\frac{\mu-\ol{\mu}}{\mu\,\ol{\mu}}\right)+2\,B_1\,\frac{\ol{x}}{\omega}\left(\mu-\ol{\mu}\right), \\
&\frac{1}{\ol{\mu}}\sqrt{R(\ol{\mu})}-\frac{1}{\mu}\sqrt{R(\mu)}
=\mathcal{I}_0\left(\frac{\mu+\ol{\mu}}{\ol{\mu}\,\mu}\right)+2\,B_1\,\ol{x}-2\,C_1\,\frac{\omega}{\ol{x}}\ ,\label{eq2:Dubrovin eq. g=1}
\end{align}\ese
where we can use that $\ol{{\kp}}={\kp}$ and fixed $\kp$ at an initial point.
Solving equation~(\ref{eq1:Dubrovin eq. g=1}) for $\ol{x}$ and inserting into~(\ref{eq2:Dubrovin eq. g=1}) we obtain the actual discrete Dubrovin equation of second-degree, namely
\be
\begin{split}
&{\left(\mu-\ol{\mu}-\omega\right)}\left(\frac{\sqrt{R(\ol{\mu})}-\mathcal{I}_0}{\ol{\mu}} \right)^2-{\left(\mu-\ol{\mu}+\omega\right)}\left(\frac{\sqrt{R({\mu})}+\mathcal{I}_0}{{\mu}} \right)^2\\
&=\frac{2\,\omega}{\mu\,\ol{\mu}}\left(\sqrt{R(\ol{\mu})}-\mathcal{I}_0\right)\left(\sqrt{R({\mu})}+\mathcal{I}_0\right)-4\,B_1\,C_1\,(\mu-\ol{\mu})^2\ ,
\end{split}
\label{Dubrovin eq. g=1}
\ee
which is different from the first-degree Dubrovin equation given in~\cite{nijhoff2000discrete} for a different discrete map.
The Casimirs $B_1$ and $C_1$ are both equal to $2\,\epsilon$.

In the case of genus $g=2$ the discrete Dubrovin equations are given by the following set of first-order difference equations
\bse \label{Dubrovin eq. g=2}
\begin{align}\nonumber
&\frac{(1/\mu_1)\sqrt{R(\mu_1)}-(1/\mu_2)\sqrt{R(\mu_2)}}{\mu_1-\mu_2}
-\frac{(1/\ol{\mu}_1)\sqrt{R(\ol{\mu}_1)}-(1/\ol{\mu}_2)\sqrt{R(\ol{\mu}_2)}}{\ol{\mu}_1-\ol{\mu}_2}\\
&= \mathcal{I}_0\left(\frac{1}{{\mu}_1{\mu}_2}+\frac{1}{\ol{\mu}_1\ol{\mu}_2}\right)
-2\,B_2\,\ol{x}+2\,C_2\,\frac{\omega}{\ol{x}}\ , \\[0.1cm] \nonumber
&\frac{(1/\mu_1)\sqrt{R(\mu_1)}-(1/\mu_2)\sqrt{R(\mu_2)}}{\mu_1-\mu_2}
+\frac{(1/\ol{\mu}_1)\sqrt{R(\ol{\mu}_1)}-(1/\ol{\mu}_2)\sqrt{R(\ol{\mu}_2)}}{\ol{\mu}_1-\ol{\mu}_2}\\
&= \mathcal{I}_0\left(\frac{1}{{\mu}_1{\mu}_2}-\frac{1}{\ol{\mu}_1\ol{\mu}_2}\right)
+2\,B_2\,\frac{\ol{x}}{\omega}\left(\mu_1-\ol{\mu}_1+\mu_2-\ol{\mu}_2\right), \\[0.1cm] \nonumber
&\frac{(\mu_2/\mu_1)\sqrt{R(\mu_1)}-(\mu_1/\mu_2)\sqrt{R(\mu_2)}}{\mu_1-\mu_2}+\frac{(\ol{\mu}_2/\ol{\mu}_1)\sqrt{R(\ol{\mu}_1)}-(\ol{\mu}_1/\ol{\mu}_2)\sqrt{R(\ol{\mu}_2)}}{\ol{\mu}_1-\ol{\mu}_2}\\
&= \mathcal{I}_0\left(\frac{\mu_1+\mu_2}{\mu_1\mu_2}-\frac{\ol{\mu}_1+\ol{\mu}_2}{\ol{\mu}_1\ol{\mu}_2}\right)
+2\,B_2\,\frac{\ol{x}}{\omega}\left(\mu_1\mu_2-\ol{\mu}_1\ol{\mu}_2\right) ,
\end{align}\ese
where the Casimirs $B_2$ and $C_2$ are both equal to $3\,\epsilon$.
As before eliminating $\ol{x}$ from (\ref{Dubrovin eq. g=2}), we get a coupled first-order difference equations of degree two for the variables $\mu_1, \mu_2$. In general, this system is difficult to solve directly. It is conceivable that the system can be solved using the techniques of Abel's paper~\cite{abel2012sur1}, but we will not pursue this line of investigation here.
Instead, we consider the integration by means of continuous-time interpolating flow. Note that equations~(\ref{Dubrovin eq. g=1}) and (\ref{Dubrovin eq. g=2}) describe the discrete evolution of the separation variables.


\subsection{Interpolating flow}\label{sec. of intrp. flow}

We will now establish the continuous evolution generated by the invariants. First, we need to establish the Poisson brackets between canonical separation variables. Using~(\ref{bracket of A1B2}) and (\ref{bracket of B1D2}) we obtain
\be
\{A(\ld_1)+D(\ld_1), B(\ld_2)\}=\frac{B(\ld_1)\left[A(\ld_2)-D(\ld_2)\right]-B(\ld_2)\left[A(\ld_1)-D(\ld_1)\right]}{\ld_1-\ld_2}\label{bracket of (A1+D1),B2} \ .
\ee
This yields in the limit $\ld_2\longrightarrow \ld_1$ the relation
\be
\{A(\ld)+D(\ld), B(\ld)\}=\left[A(\ld)-D(\ld)\right]B'(\ld)-B(\ld)\left[A(\ld)-D(\ld)\right]'\ ,
\label{eq: Poisson of tr(T) with B}
\ee
where the prime denotes the differentiation with respect to $\ld$.
From the auxiliary spectrum (\ref{eq:auxiliary spectrum B}), using $B(\mu_j)=0$ and the relation
\begin{align*}
\{A+D, B\}(\mu_j)=(A-D)(\mu_j)\,B'(\mu_j)\ ,
\end{align*}
we obtain
\be
\{(A+D)(\mu_j), \mu_j\}=(A-D)(\mu_j)\ ,
\ee
whereas for $\mu_i\neq\mu_j$ we have $$\{\mu_i, (A+D)(\mu_j)\}=0\ .$$
A set of canonical separation variables is introduced by taking the $\{\mu_j\}$ as the position variables, and $\{\nu_j\}$ as momenta variables defined by
\be
2\,\textrm{cosh}(\nu_j)\equiv\frac{{\rm tr}\,T(\mu_j)}{\sqrt{\det T(\mu_j)}}\quad \Rightarrow\quad \nu_j=\hf\,\textrm{log}\left(\frac{A(\mu_j)}{D(\mu_j)}\right)\ ,
\label{eq:nu and cosh def}
\ee
leading to the Poisson brackets
\begin{equation}\label{eq: canon. with mu,nu}
\{ \mu_{i},\mu_{j}\}=\{\nu_{i},\nu_{j}\}=0, \qquad \{\mu_{i}, \nu_{j}\}=\dd_{ij} ,\qquad j=1,\dots,g\ .
\end{equation}

In order to construct an interpolating flow for the map~(\ref{map T (SoV)}), we consider the continuous-time evolution for the auxiliary spectrum using $\textrm{tr}\,T(\ld)$ as the generating Hamiltonian of the flow\footnote[12]{In choosing a single parameter-flow generated by $\textrm{tr}\,{T}({\ld})$ we anticipate matching this flow to the map in the case of $g=2$. For higher genus we would need a multiparameter combination of traces. We will not consider the latter extension in this paper.} (a similar time evolution was considered in~\cite{cao2012finite} by Cao \& Xu), i.e.
\be
\dot{\mu}_j=\frac{\partial\mu_j}{\partial t_{\ld}}=\{\mu_j, \textrm{tr}\,{T}({\ld})\}=\frac{\sqrt{R(\mu_j)}\,B(\ld)}{(\ld-\mu_j)\,B_g\prod_{i\neq j}(\mu_j -\mu_i)}\ ,
\label{eq:Mu dot(g=P-1)}
\ee
which are the Dubrovin equations for our case. The value of $\ld$ will be fixed later.
From the definition~(\ref{eq:nu and cosh def}) and equation~(\ref{eq:Mu dot(g=P-1)}), the second companion equation for the conjugate variable ${\nu}_j$ is given by
\be
\dot{\nu}_j =\frac{\sqrt{\det T(\mu_j)}\,B(\ld)}{(\ld-\mu_j)\,B'(\mu_j)}\,\frac{\partial}{\partial \mu_j}\left(\frac{{\rm tr}\,T(\mu_j)}{\sqrt{\det T(\mu_j)}}\right).
\label{eq:Nu dot(g=P-1)}\ee
The coupled equations~(\ref{eq:Mu dot(g=P-1)})~and~(\ref{eq:Nu dot(g=P-1)}) derive from the Hamiltonian of the form
\be
H_{\ld}(\mu_1, \dots, \mu_{g}; \nu_1, \dots, \nu_{g})=\sum_{j=1}^{g}\frac{2\,\textrm{cosh}(\nu_j)\sqrt{\det T(\mu_j)}\,B(\ld)-B(\ld)\,\textrm{tr}\,{T}(\mu_j)}{(\ld-\mu_j)\,B'(\mu_j)}\ ,
\label{eq:Hamiltonian_lambda(SoV:g=P-1)}
\ee
which is different from Toda's equation, cf.~\cite{toda2012theory}, where he gave the potential term only.
Using the Lagrange interpolation formula
\begin{equation}
\frac{\ld^k}{B(\ld)}+\sum_{i=1}^g\frac{\mu_i^k}{B'(\mu_i)(\mu_i-\ld)}=\begin{cases}
	\,0 \,\quad\  \ , \quad 0\leq k\leq g-1 \\
	{1}/{B_g}\ , \quad k=g
\end{cases}\ ,
\label{Lagrange interpol.}
\end{equation}
we can integrate the Dubrovin equations (\ref{eq:Mu dot(g=P-1)}), leading to the following Jacobi inversion problem
\be
\ld^k(t-t_0)=\sum_{j=1}^{g} \int_{\mu_j(t_0)}^{{\mu}_j(t)}\, \frac{\mu^k}{\sqrt{R(\mu)}}\, d\mu,\qquad k=0,\dots, g-1\ ,
\label{t:int.for g=P-1:SoV}
\ee
on the hyperelliptic Riemann surface of genus $g=P-1$
\be
\Gamma: \quad  \upeta^2-\tr T(\ld)\,\upeta+\det T(\ld)=0\ .
\label{eq:spectral curve in terms of tr, det}
\ee
We note that the time variable $t=t(\ld)$ depends on the parameter value $\ld$.
We also note that for different values of $\ld$ the compatibility flows commute as a consequence of equation~(\ref{eq: involution property}).


\subsection{Generating functions structures and action-angle variables}\label{sec. of Generating functions}

We can now introduce the canonical transformation to action-angle variables in terms of the spectral variables $\mu_1, \dots, \mu_g, \nu_1, \dots, \nu_g$. We look for a canonical transformation given by a generating function $G(\mu_1, \dots, \mu_g; \mathcal{P}_1, \dots, \mathcal{P}_g)$ in which the action variables $\mathcal{P}_1, \dots, \mathcal{P}_g$ are the invariants appearing in~(\ref{trace for KdV type (SoV)}); specifically $\mathcal{P}_j=\mathcal{I}_{j-1}$ where $j=1,\dots,g$. This generating function is obtained from
\be
G=\sum_{j=1}^{g} \int_{\mu_j^0}^{{\mu}_j} \nu_j(\mu;\mathcal{P}_1, \dots, \mathcal{P}_g) \, d\mu=\sum_{j=1}^{g} \int_{\mu_j^0}^{{\mu}_j}\textrm{arccosh}\left(\hf\frac{\tr T(\mu)}{\sqrt{\det T(\mu)}}\right) d\mu\ ,
\label{eq:G function (general)}
\ee
up to an arbitrary function of the invariants. Note that the integral~(\ref{eq:G function (general)}) depends on $\mu$'s and $\mathcal{P}$'s only.
The angle variables are obtained from
\be
\mathcal{Q}_{k}=\frac{\partial G}{\partial \mathcal{P}_{k}}=\sum_{j=1}^{g} \int_{\mu_j(t_0)}^{{\mu}_j(t)} \frac{\mu^{k-1}}{\sqrt{R(\mu)}}\, d\mu, \qquad k=1, \dots, g\ .
\ee

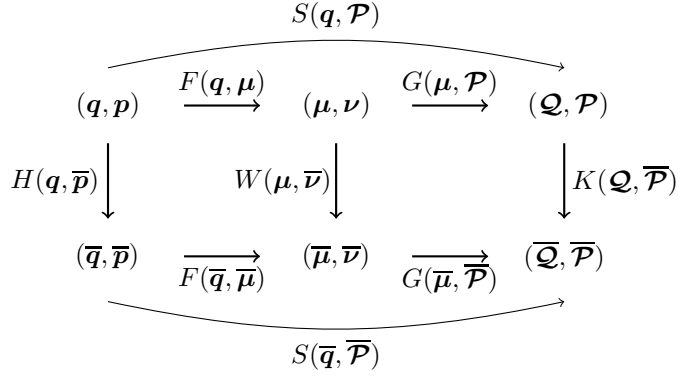
\begin{figure}[ht]
\centering
\begin{tikzpicture}[every node/.style={minimum size=1cm}]
\draw[thick,->] (1,2) -- (2,2) node[right=11] {$(\boldsymbol{\mathcal{Q}}, \boldsymbol{\mathcal{P}})$};
\draw[thick,->] (-2,2) -- (-1,2) node[right=12] {$(\boldsymbol{\mu}, \boldsymbol{\nu})$};
\draw[thick,->] (1,0) -- (2,0) node[below left=-0.3] {};
\draw[thick,<-] (-1,0) -- (-2,0) node[left=13] {$(\ol{\boldsymbol{q}}, \ol{\boldsymbol{p}})$};
\draw[thick,->] (3,1.5) -- (3,0.5) node[below=-0.3] {$(\ol{\boldsymbol{\mathcal{Q}}}, \ol{\boldsymbol{\mathcal{P}}})$};
\draw[thick,->] (0,1.5) -- (0,0.5) node[below=-0.3] {$(\ol{\boldsymbol{\mu}}, \ol{\boldsymbol{\nu}})$};
\draw[thick,<-] (-3,0.5) -- (-3,1.5) node[above=-0.3] {$(\boldsymbol{q}, \boldsymbol{p})$};
\draw[<-] (3,2.5) to [bend right=12] (-3,2.5) node[below=-0.3] {};
\draw[->] (-3,-0.6) to [bend right=12] (3,-0.6);
\node at (-1.5,2.3) {$F(\boldsymbol{q}, \boldsymbol{\mu})$};
\node at (-1.5,-0.3) {$F(\ol{\boldsymbol{q}}, \ol{\boldsymbol{\mu}})$};
\node at (1.5,2.3) {$G(\boldsymbol{\mu}, \boldsymbol{\mathcal{P}})$};
\node at (1.5,-0.3) {$G(\ol{\boldsymbol{\mu}}, \ol{\boldsymbol{\mathcal{P}}})$};
\node at (-3.7,1) {$H(\boldsymbol{q}, \ol{\boldsymbol{p}})$};
\node at (-0.7,1) {$W({\boldsymbol{\mu}}, \ol{\boldsymbol{\nu}})$};
\node at (3.8,1) {$K(\boldsymbol{\mathcal{Q}}, \ol{\boldsymbol{\mathcal{P}}})$};
\node at (0,3.2) {$S(\boldsymbol{q}, \boldsymbol{\mathcal{P}})$};
\node at (0,-1.3) {$S(\ol{\boldsymbol{q}}, \ol{\boldsymbol{\mathcal{P}}})$};
\end{tikzpicture}
\caption{Commuting diagram of canonical transformations}
\label{commuting diagram of integrable canonical}
   \end{figure}

As a consequence we have that the transformation to action-angle variables on the one hand and the mapping on the other hand form a ladder of commuting canonical transformations, as shown in figure~\ref{commuting diagram of integrable canonical}.
In fact, we can think of the transformation between the original variables and action-angle variables as the result of the composition of two generating functions which factorize the canonical transformation that leads to the MH. To our knowledge, there is no existing general theory of composition of generating functions, and we basically give a description by the following commuting diagram for what we mean by that.
   
In this diagram, the $S$ denotes the generating function of the canonical transformation from the original variables to the action-angle variables, the $F$ denotes the generating function of the canonical transformation from the original variables to the separating variables, and the $G$ denotes the generating function of the canonical transformation from the separating variables to the action-angle variables.
The $H$ is the action functional describing the canonical transformation which is the discrete mapping (Hamiltonian). The canonical transformation, with generating function $W$, realising the dynamical mapping in terms of the separation variables, is obviously, by construction, an integrable map itself. Unfortunately, however, it does not seem easy in general to obtain an explicit expression for the generating function $W$, since this requires the elimination of the invariants entering as coefficients of the spectral curve. The $K$ is the action functional describing the canonical transformation in terms of the action-angle variables.

As an explicit example of the above diagram, let us consider the KdV map example~\ref{sec. of KdV example(g=1)}. In this example the generating function $S$ satisfying the defining equations (\ref{eq:action-angle defining eqs(q,p)KdV,g=1)}) is given by (\ref{eq:S function (KdV g=1)}).
The generating functions $F$ and $G$ which are, respectively, parametrized as
\bse
\begin{align}
{p}&=\frac{\partial{F}}{\partial{{q}}}, \qquad\, {\nu}=-\frac{\partial{F}}{\partial{\mu}}\ ,\\
{\nu}&=\frac{\partial{G}}{\partial{{\mu}}}, \qquad {\mathcal{Q}}=\frac{\partial{G}}{\partial{\mathcal{P}}}\ ,
\end{align}\label{defining eqs:for F,G}
\ese
are given by
\begin{align}
F(q,\mu)=&\,\left(\omega -\mu\right)\textrm{log}\,(\epsilon+q)+\mu\,\textrm{log}\,(\epsilon-q)-\epsilon\,q\ , 	\label{eq:F function(KdV)g=1} \\
G(\mu, \mathcal{P})=&\int_{\mu^0}^{\mu}
\textrm{log}\,{\left(\frac{\sqrt{4\,\epsilon^{2}\mu'+\mathcal{P}+\alpha^{2}}
+\sqrt{4\,\mu'\left(\mu'+\alpha\right)+\mathcal{P}+\alpha^{2}}}{\sqrt{4\,\epsilon^{2}\mu'+\mathcal{P}+\alpha^{2}}
-\sqrt{4\,\mu'\left(\mu'+\alpha\right)+\mathcal{P}+\alpha^{2}}}\right)}\,d\mu'\ ,
	\label{eq:G function(KdV)g=1}
\end{align}
up to an arbitrary function of the invariant, and in which $\alpha\equiv \epsilon^2-\omega$. In equations~(\ref{defining eqs:for F,G}) $\mu$ and $\nu$ are given  by
\begin{align}
\mu=&\,\frac{1}{2\,\epsilon}\left(\epsilon\,\omega-\omega\,q-{\epsilon}^{3}-{\epsilon}^{2}p+\epsilon\,{q}^{2}+p\,{q}^{2}\right),\\
\nu=&\,\textrm{log}\left(\frac{\epsilon+q}{\epsilon-q}\right),
\label{eq: mu and nu KdV g=1}
\end{align}
where $\mu$ is the auxiliary spectrum and $\nu$ defined in~(\ref{eq:nu and cosh def}).
The generating function $H$ satisfying the relations (\ref{discrete Hamilton's equations KdV g=1}) is given by (\ref{Hamiltonian eq. for KdV case P=2)}),
whereas the generating function $K$ satisfying the relations
$$\ol{\mathcal{P}}-{\mathcal{P}}=-\frac{\partial{K}}{\partial{{\mathcal{Q}}}}, \qquad \ol{{\mathcal{Q}}}-{\mathcal{Q}}=\frac{\partial{K}}{\partial{\mathcal{P}}}\ ,$$
is given by~(\ref{eq: closed form-(KdV g=1)}),
which is in fact the MH of the map~(\ref{KdV map in q,p(g=1)}).
We note that the generating function $K$ is presented in terms of the variable $\mathcal{P}$ only since $\mathcal{P}$ is an invariant, i.e.\ $\ol{\mathcal{P}}=\mathcal{P}$. The schemes that we set up in section~\ref{sec. of SoV for KdV} will be illustrated for genus-two in the next section.


\section{The modified Hamiltonian of two-degrees-of-freedom}\label{sec. of 2-degrees Hamiltonian}
\setcounter{equation}{0}

As discussed in section~\ref{sec. of KdV and mKdV mappings}, the lattice KdV equation~(\ref{KdV eq}) leads to the integrable mapping~(\ref{mapping of KdV (P=3)in X and Y}) when considering a staircase with period $P=3$. For convenience, this mapping is introduced below:
\begin{equation}
\begin{split}
		&\overline{p}_1=p_1+\frac{\epsilon\, \delta}{\epsilon-q_1}-\frac{\epsilon\, \delta}{\epsilon+q_1+q_2} , \quad
		\overline{q}_1=q_1-\frac{\epsilon\, \delta}{\epsilon-\overline{p}_1}+\frac{\epsilon\, \delta}{\epsilon+\overline{p}_1-\overline{p}_2}\ , \\
		&\overline{p}_2=p_2+\frac{\epsilon\, \delta}{\epsilon-q_2}-\frac{\epsilon\, \delta}{\epsilon+q_1+q_2} , \quad
		\overline{q}_2=q_2+\frac{\epsilon\, \delta}{\epsilon+\overline{p}_2}-\frac{\epsilon\, \delta}{\epsilon+\overline{p}_1-\overline{p}_2}\ ,
		\end{split}
		\label{KdV mapping (P=3)} 
\end{equation}
in which we identify $$X_1:= p_1,\qquad X_2:= p_2-p_1,\qquad Y_1:= q_1,\qquad Y_2:= q_2\ .$$
A further discussion of this mapping can be found in~\ref{appendix:b}.
In terms of the conjugate variables $q_1, q_2, p_1, p_2$, we have the standard symplectic structure
\begin{equation}
\Omega=dq_1 \wedge dp_1+dq_2 \wedge dp_2\ ,
\label{symplectic structure KdV mapping (P=3)}
\end{equation}
leading to the standard Poisson brackets
$$\{q_1,q_{2}\}=\{p_1,p_{2}\}=\{q_1,p_{2}\}=\{q_2,p_{1}\}=0,\quad \{q_1,p_{1}\}=\{q_2,p_{2}\}=1\ .$$
This structure is preserved by the map (\ref{KdV mapping (P=3)}), and the mapping is in fact a canonical transformation, i.e.\ $\overline{\Omega}=\Omega$.
The Hamiltonian,
\begin{equation}
\begin{split}
    {H}(q_1,q_2,\overline{p}_1,\overline{p}_2)=&\,\epsilon \,\delta\, \textrm{log}\,(\epsilon+\overline{p}_1-\overline{p}_2)(\epsilon-\overline{p}_1)(\epsilon+\overline{p}_2) \\
    &+\epsilon \,\delta\, \textrm{log}\,(\epsilon+q_1+q_2)(\epsilon-q_1)(\epsilon-q_2)\ ,
    \end{split}
	\label{Hamiltonian eq. for KdV P=3}
\end{equation}
acts as the generating function for the mapping (\ref{KdV mapping (P=3)}), i.e. one has the discrete-time Hamilton equations
\begin{equation*}
\begin{split}
	&\overline{p}_1-p_1=-\frac{\partial H}{\partial q_1}, \qquad   \ \overline{q}_1-q_1=\frac{\partial H}{\partial \overline{p}_1}\ , \\
	&\overline{p}_2-p_2=-\frac{\partial H}{\partial q_2}, \qquad   \ \overline{q}_2-q_2=\frac{\partial H}{\partial \overline{p}_2}\ ,
	\end{split}
	\label{Hamilton's equations for KdV (P=3)}
\end{equation*}
leading to the equations of the mapping (\ref{KdV mapping (P=3)}).
As noted around~(\ref{eq:invariants(X,Y):KdV g=2}), the mapping~(\ref{KdV mapping (P=3)}) conserves the invariants
\bse
\begin{align} \nonumber
\mathcal{P}_1=&\,\frac{1}{8}\,{\delta}^{2}\,{\epsilon}^{2} \left\lbrack {x_1} \left(3\,{y_0}+3\,{y_1}-{y_2}\right)+{x_2}\left(3\,{y_0}+3\,{y_2}-{y_1} \right)+{x_3} \left(3\,{y_1}+3\,{y_2}-{y_0} \right)  \right\rbrack\\ \nonumber
&+\frac{1}{2}\,\delta\,\epsilon\left\lbrack {x_1}\,{x_2}\,{y_0}\left({y_2}-{y_1} \right)+{x_1}\,{x_3}\,{y_1}\left( {y_0}-{y_2}\right)+{x_2}\,{x_3}\,{y_2} \left( {y_1}-{y_0} \right)\right\rbrack\\
&+{x_1}\,{x_2}\,{x_3}\,{y_0}\,{y_1}\,{y_2}\ ,\label{eq:1st invariant KdV(g=2)} \\[0.1cm] \nonumber
\mathcal{P}_2=&\,\delta\,\epsilon \left\lbrack {x_1} \left( {y_0}-{y_1} \right)+{x_2}\left({y_2}-{y_0} \right) +{x_3} \left({y_1}-{y_2} \right) \right\rbrack\\
&+{x_1}\,{x_2}\,{y_0}\left( {y_1}+{y_2} \right)+{x_1}\,{x_3}\,{y_1} \left({y_0}+{y_2}\right)+{x_2}\,{x_3}\,{y_2}\left( {y_0}+{y_1}\right),
\end{align}
   \label{eq:invariants(q,p):KdV g=2}
\ese
in which we have used the abbreviations
$$x_1\equiv\epsilon-p_1,~\, x_2\equiv\epsilon+p_2,~\, x_3\equiv\epsilon+p_1-p_2,~\, y_0\equiv\epsilon+q_1+q_2,~\, y_1\equiv\epsilon-q_1,~\, y_2\equiv\epsilon-q_2\ .$$
The canonical structure allows us to show the integrability property that the two invariants are in involution with each other, with respect to the canonical Poisson bracket:
$\{\mathcal{P}_1, \mathcal{P}_2\}=0$.
The invariance and involutivity of these can be shown by direct calculations. The calculations involved are large, but have been verified by MAPLE.
The invariants $\mathcal{P}_1$ and $\mathcal{P}_2$ will thus generate two commuting continuous flows to the mapping~(\ref{KdV mapping (P=3)}).

Now, by applying the BCH formula~(\ref{eq:BCH for MH (general)}) to the Hamiltonian (\ref{Hamiltonian eq. for KdV P=3}), where we consider $\delta$ as the step size of the mapping (\ref{KdV mapping (P=3)}), one obtains the following expression for the modified Hamiltonian:
\begin{equation}
	\begin{split}
		{H}^{\ast}=&\,\delta\,\epsilon\,\textrm{log}\,(\epsilon+p_1-p_2)(\epsilon-p_1)(\epsilon+p_2)+\delta\,\epsilon\,\textrm{log}\,(\epsilon+q_1+q_2)(\epsilon-q_1)(\epsilon-q_2)\\
		&-\frac{\delta^2\,\epsilon^2}{2}\left[ \left(\frac{1}{\epsilon+p_1-p_2}-\frac{1}{\epsilon-p_1}\right)\left(\frac{1}{\epsilon+q_1+q_2}-\frac{1}{\epsilon-q_1}\right) \right. \\
		&\qquad~ \left. -\left(\frac{1}{\epsilon+p_1-p_2}-\frac{1}{\epsilon+p_2}\right)\left(\frac{1}{\epsilon+q_1+q_2}-\frac{1}{\epsilon-q_2}\right)\right] 
		+O(\delta^3)\ .
	\end{split}
	\label{eq:MH: KdV(g=2):BCHseries}
\end{equation}
We can also express the MH in terms of the invariants:
\be
{H}^{\ast}=\delta\,\epsilon\,\textrm{log}\left(\mathcal{P}_1 \right)+O(\delta^3)\ .
\label{eq:MH in terms invariants(KdV, g=2)}
\ee
After inserting the invariant~$\mathcal{P}_1$ and expanding equation~(\ref{eq:MH in terms invariants(KdV, g=2)}) is seen to agree with~(\ref{eq:MH: KdV(g=2):BCHseries}) up to order $\delta^{3}$.
Observe that in equation~(\ref{eq:MH in terms invariants(KdV, g=2)}), the term corresponding to $\delta^{2}$ is zero. In fact, in the light of our present investigation we expect that for Hamiltonian system of KdV models, all terms of the interpolating Hamiltonian which correspond to $\delta^{2n}$ where $n\in\mathbb N$, are equal to zero.

Turning to the action-angle variables technique, accordance used Hamilton-Jacobi theory we need separation variables to obtain a closed-form expression for the MH (\ref{eq:MH: KdV(g=2):BCHseries}). This is advantageous in this case because in this manner we can take care of the separation of variables.
$\mathcal{P}_1$ and $\mathcal{P}_2$ can be defined as the new momenta and their canonically conjugated variables as the new coordinates $\mathcal{Q}_1, \mathcal{Q}_2$,
thus, we have a canonical transformation
\be({\mu}_1, \mu_2,{\nu}_1, \nu_2) \longrightarrow ({\mathcal{Q}}_1, \mathcal{Q}_2,{\mathcal{P}}_1, \mathcal{P}_2)\label{canonical transformation for g=2} \ee
which can be given by means of a generating function, $G(\mu_1, \mu_2, \mathcal{P}_1, \mathcal{P}_2)$ as
\be
{\nu}_1=\frac{\partial{G}}{\partial{{\mu}_1}}, \quad {\nu}_2=\frac{\partial{G}}{\partial{{\mu}_2}}, \quad \mathcal{Q}_1=\frac{\partial{G}}{\partial{\mathcal{P}_1}}, \quad \mathcal{Q}_2=\frac{\partial{G}}{\partial{\mathcal{P}_2}}\ ,
\label{eq:defining eqs(mu,nu)KdV(g=2)}\ee
with
\be\label{eq:K_lambda(KdV, g=2)}
	K_{\ld}=H_{\ld}+\frac{\partial G}{\partial t_{\ld}}\ ,
\ee
being the transformed Hamiltonian.
Thus, we get $G$ up to a function of the invariants
\be
G=\int_{\mu_1^0}^{{\mu}_1}\textrm{arccosh}\left(\hf\frac{\tr T(\mu)}{\sqrt{\det T(\mu)}}\right) d\mu+\int_{\mu_2^0}^{{\mu}_2}\textrm{arccosh}\left(\hf\frac{\tr T(\mu)}{\sqrt{\det T(\mu)}}\right) d\mu\ ,
\label{eq: G function g=2}
\ee
and consequently we have
\bse\begin{align}
&\mathcal{Q}_1=\int_{\mu_1^0}^{{\mu}_1} \frac{{1}}{\sqrt{R(\mu)}}\, d\mu+\int_{\mu_2^0}^{{\mu}_2} \frac{1}{\sqrt{R(\mu)}}\, d\mu\ ,\\
&\mathcal{Q}_2=\int_{\mu_1^0}^{{\mu}_1} \frac{\mu}{\sqrt{R(\mu)}}\, d\mu+\int_{\mu_2^0}^{{\mu}_2} \frac{\mu}{\sqrt{R(\mu)}}\, d\mu\ .
\end{align}\ese
We note that in the case of genus two, the discriminant of the curve takes on the form:
\begin{equation}
	\begin{split}
R(\ld)=&\,36\,\epsilon^2\,\ld^5+\left(4\,\mathcal{P}_2-3\,\omega^2-54\,\epsilon^2\,\omega+81\,\epsilon^4 \right)\ld^4 \\
       &+2\left(9\,\epsilon^2\,\mathcal{P}_2-3\,\omega\,\mathcal{P}_2+2\,\mathcal{P}_1+2\,\omega^3\right)\ld^3\\
       &+\left(18\,\epsilon^2\,\mathcal{P}_1-6\,\omega\,\mathcal{P}_1+{\mathcal{P}_2}^2\right)\ld^2
       +2\,\mathcal{P}_1\,\mathcal{P}_2\,\ld+{\mathcal{P}_1}^2\ .
	\end{split}
	\label{eq:disc(KdV)g=2}
\end{equation}

As discussed in section~\ref{sec. of intrp. flow}, a continuous-time evolution for the auxiliary spectrum can be introduced by considering $\textrm{tr}\,T(\ld)$ as the generating Hamiltonian of the flow. In the two-degrees-of-freedom case, the trace of the monodromy matrix is given by
\be
\textrm{tr}\,T(\ld) = 2\,\lambda^3+3\left(3\,\epsilon^2-\omega\right)\lambda^2+\mathcal{P}_2\,\lambda+\mathcal{P}_1\ .
\label{eq:trace KdV g=2}
\ee
Using the trace~(\ref{eq:trace KdV g=2}) as the generating Hamiltonian of the flow, one can obtain the following set of equations
\bse\label{eq:Hamilton Mu,Nu (g=2)}
\begin{align}
&\dot{\nu}_1 = \frac{(\ld-\mu_2)\sqrt{\det T(\mu_1)}}{\mu_1-\mu_2}\,\frac{\partial}{\partial \mu_1}\left(\frac{{\rm tr}\,T(\mu_1)}{\sqrt{\det T(\mu_1)}}\right),\quad\,
\dot{\mu}_1 = \frac{(\ld-\mu_2)\sqrt{R(\mu_1)}}{\mu_1-\mu_2}\ , \\
&\dot{\nu}_2 = \frac{(\ld-\mu_1)\sqrt{\det T(\mu_2)}}{\mu_2-\mu_1}\,\frac{\partial}{\partial \mu_2}\left(\frac{{\rm tr}\,T(\mu_2)}{\sqrt{\det T(\mu_2)}}\right),\quad\,
\dot{\mu}_2 = \frac{(\ld-\mu_1)\sqrt{R(\mu_2)}}{\mu_2-\mu_1}\ ,
\end{align}\ese
which derive from the Hamiltonian of the form
\be
\begin{split}
H_{\ld}(\mu_1, \mu_2, \nu_1, \nu_2)=\sum_{j=1,2}\prod_{i\neq j}(\ld-\mu_i)\left(\frac{2\,\textrm{cosh}(\nu_j)\sqrt{\det T(\mu_j)}-\tr T(\mu_j)}{\prod_{i\neq j}(\mu_j -\mu_i)} \right).
\end{split}
\label{Hamiltonian SoV:g=2}
\ee
Note that the dot $\cdot$ in~(\ref{eq:Hamilton Mu,Nu (g=2)}) denotes the differentiation with respect to the continuous-time flow variable $t_\ld$.
Using (\ref{Lagrange interpol.}) we can integrate (\ref{eq:Hamilton Mu,Nu (g=2)}) to obtain
\bse\label{eq:t-integral(KdV, g=2, SoV)}\begin{align}
(t-t_0)&=\int_{\mu_1(t_0)}^{{\mu}_1(t)}\, \frac{1}{\sqrt{R(\mu)}}\, d\mu+\int_{\mu_2(t_0)}^{{\mu}_2(t)}\, \frac{1}{\sqrt{R(\mu)}}\, d\mu\ , \\
\ld(t-t_0)&=\int_{\mu_1(t_0)}^{{\mu}_1(t)}\, \frac{\mu}{\sqrt{R(\mu)}}\, d\mu+\int_{\mu_2(t_0)}^{{\mu}_2(t)}\, \frac{\mu}{\sqrt{R(\mu)}}\, d\mu\ .
\end{align}\ese
In fact, the modified Hamiltonian ${H}^{\ast}$ coincides with the canonical transformed Hamiltonian $K_\ld$ obtained by applying the canonical transformation (\ref{canonical transformation for g=2}), viewed as a function of $\mathcal{Q}_1, \mathcal{Q}_2, \mathcal{P}_1, \mathcal{P}_2$. Hamilton's equations in the new coordinates are given by
\bse\label{eq:Hamilton's eqs (P1,P2,Q1,Q2}
\begin{align}
	&\dot{\mathcal{P}}_1=-\frac{\partial {H}^{\ast}}{\partial \mathcal{Q}_1}=0, \qquad \dot{\mathcal{Q}}_1=\frac{\partial {H}^{\ast}}{\partial \mathcal{P}_1}=\Nu_1\ , \\
	&\dot{\mathcal{P}}_2=-\frac{\partial {H}^{\ast}}{\partial \mathcal{Q}_2}=0, \qquad \dot{\mathcal{Q}}_2=\frac{\partial {H}^{\ast}}{\partial \mathcal{P}_2}=\Nu_2\ ,
\end{align}\ese
which tell us that ${H}^{\ast}$ depends on $\mathcal{P}_1, \mathcal{P}_2$ only.

Using Abel's theorem, cf.\ \ref{appendix:e} for more discussion, we introduce the frequencies $\Nu_1$, $\Nu_2$ (not to be confused with a canonical variables $\nu_1,\nu_2$) as the discrete time-one step as the following
\bse\label{eqs:freq. by Abel,g=2}
\begin{align}
\Nu_1=\int_{\mu_1}^{\ol{\mu}_1}\, \frac{1}{\sqrt{R(\mu)}}\, d\mu+\int_{\mu_2}^{\ol{\mu}_2}\, \frac{1}{\sqrt{R(\mu)}}\, d\mu=-\int_{\infty}^{(\omega,\upeta(\omega))}\frac{1}{\sqrt{R(\mu)}}\, d\mu\ ,\\
\Nu_2=\int_{\mu_1}^{\ol{\mu}_1}\, \frac{\mu}{\sqrt{R(\mu)}}\, d\mu+\int_{\mu_2}^{\ol{\mu}_2}\, \frac{\mu}{\sqrt{R(\mu)}}\, d\mu=-\int_{\infty}^{(\omega,\upeta(\omega))}\frac{\mu}{\sqrt{R(\mu)}}\, d\mu\ ,
\end{align}\ese
so that $$\ol{\mathcal{Q}}_1-\mathcal{Q}_1=\Nu_1,\qquad\textrm{and}\qquad \ol{\mathcal{Q}}_2-\mathcal{Q}_2=\Nu_2\ ,$$ and in which we have chosen a specific value of $\ld$ given by
\be\label{eq:value of lamda}
\ld=\left({\int_{\infty}^{(\omega,\upeta(\omega))}\frac{\mu}{\sqrt{R(\mu)}}\, d\mu}\right)/\left({\int_{\infty}^{(\omega,\upeta(\omega))}\frac{1}{\sqrt{R(\mu)}}\, d\mu}\right) .
\ee
In fact, the integrals~(\ref{eqs:freq. by Abel,g=2})
provide us with the solution of the discrete Dubrovin equations~(\ref{Dubrovin eq. g=2}).

As the system of equations~(\ref{Dubrovin eq. g=2}) is given, we can then follow the KdV map example prescription of section~\ref{sec. of KdV example(g=1)} to compute the limits $\overline{\mu}_1$, $\overline{\mu}_2$ of the integrals~(\ref{eqs:freq. by Abel,g=2}) at some suitably chosen initial points $\mu_1$, $\mu_2$. However, the calculations involved are very large and cannot be reproduced here, but we shall just give the steps how to compute these limits.
Step one, consider the auxiliary spectrum $\mu_1$, $\mu_2$ (i.e.~the roots of polynomial~(\ref{eq:auxiliary spectrum B}) for $g=2$) which are expressed in terms of the variables $q_1, p_1, q_2, p_2$.
Step two, let $q_1= q_2=0$, and thus $\mu_1$, $\mu_2$, respectively, take the expressions
\bse  \label{eq:initial mu1,mu2}
\begin{align}
   &\mu_1 = \frac{1}{6}\left(3\,{\omega}+ 2\,\epsilon\,p_2-4\,{\epsilon}^{2}+ \sqrt{40\,\epsilon^4-4\,\epsilon^3\,p_2-8\,\epsilon^2\,p_2^2-6\,\mathcal{P}_2-3\,\omega^2}\right) , \\
   &\mu_2 = \frac{1}{6}\left(3\,{\omega}+ 2\,\epsilon\,p_2-4\,{\epsilon}^{2}- \sqrt{40\,\epsilon^4-4\,\epsilon^3\,p_2-8\,\epsilon^2\,p_2^2-6\,\mathcal{P}_2-3\,\omega^2}\right) .
\end{align} \ese
Step three, use the coupled system of equations~(\ref{eq:invariants(q,p):KdV g=2}) at $q_1= q_2=0$, i.e.
\begin{align*}
  &\mathcal{P}_1= \frac{\epsilon^2}{8}\left(4\,\mathcal{P}_2-8\,\epsilon\, p_1^{2}\,p_2+8\,\epsilon\,p_1\,p_2^{2}+15\,\omega^{2}-16\,\epsilon^{4}\right), \\
  &\mathcal{P}_2= -2\,\epsilon^{2}\left(p_1^{2}-p_1\,p_2+p_2^{2}-3\,\epsilon^{2}\right),
\end{align*}
to eliminate $p_2$ from the equations~(\ref{eq:initial mu1,mu2}).
Step four, rewrite the coupled system of equations~(\ref{eq:initial mu1,mu2}) as a system for which $\mu_1$, $\mu_2$ depend on the variables $\mathcal{P}_1$, $\mathcal{P}_2$ only.
Step five, use the system~(\ref{Dubrovin eq. g=2}) and equation~(\ref{eq:disc(KdV)g=2}) together with the initial points $\mu_1$, $\mu_2$ obtained from the latter step to compute $\overline{\mu}_1$, $\overline{\mu}_2$.

Denoting the new Hamiltonian by ${H}^{\ast}(\mathcal{P}_1, \mathcal{P}_2)$ one obtains
\be\begin{split}
{H}^{\ast}(\mathcal{P}_1, \mathcal{P}_2)
=\int^{(\mathcal{P}_{1}, \mathcal{P}_{2})}\left[\Nu_1\,d\mathcal{P}_1+\Nu_2\,d\mathcal{P}_2\right] ,
\label{eq:closed form MH (KdV g=2)}
\end{split}\ee
over any curve in the $\mathcal{P}_1, \mathcal{P}_2$-plane up to an integration constant, and where the expression between brackets is a differential $1$-form. The key reason is that on the basis of the integrals~(\ref{eqs:freq. by Abel,g=2}), using the dependence of the discriminant $R(\mu)$ on $\mathcal{P}_1$ and $\mathcal{P}_2$, one can prove that
\be
\frac{\partial \Nu_1}{\partial \mathcal{P}_2}=\frac{\partial \Nu_2}{\partial \mathcal{P}_1}\ ,
\label{eq:derivatives of Nu_1,Nu_2}
\ee
asserting that $(\Nu_1,\Nu_2)$ is a conservative vector field. Hence the integral~(\ref{eq:closed form MH (KdV g=2)}) is independent of the integration path in $\mathcal{P}_1, \mathcal{P}_2$-plane chosen and thus leads to a well-defined function of $(\mathcal{P}_1, \mathcal{P}_2)$ obeying the relations~(\ref{eq:Hamilton's eqs (P1,P2,Q1,Q2}).
Equation~(\ref{eq:closed form MH (KdV g=2)}) is a closed-form expression for the MH of the mapping (\ref{KdV mapping (P=3)}).


\section{Conclusion}\label{sec. of Conclusion}

We have presented examples of integrable numerical integration schemes arising from the reduction of nonlinear integrable lattice equations, which have closed form MHs. Furthermore, we have discussed the extension to multiple-degrees-of-freedom. It is generally understood that common numerical methods applied to Hamiltonians of the Newtonian form $H=p^2/2+V(q)$, where $\partial\,V(q)/\partial\,q$ is nonlinear, have an expansion for the MH which does not converge. However, in the light of the present study, there exist special symplectic integrators that do allow convergent expansion for the MHs for some non-Newtonian Hamiltonian systems, which are associated with discrete integrable dynamics. The corresponding Hamiltonian systems are associated with the interpolating flow of these integrable mappings. We have presented an example of the Hamiltonian system for the simplest case (i.e.\ one-degree-of-freedom) which arises from the reduction of the lattice KdV equation. In this case, we obtained a closed-form expression in terms of an elliptic integral and compare it with the expansion that we get from Yoshida's approach.
 
In the present paper, we broadened the perspective by looking at two parallel extensions. One was the extension to an implicit scheme where the MKdV case is key example. In particular, the example of the MKdV map exhibits an implicit dependence on the time step which could be of relevance to certain implicit schemes in numerical analysis. Another extension was the extension to multiple-degrees-of-freedom. In the latter case, we are dealing with the more complicated situation where the underlying spectral representation of the underlying of integrable systems is associated with higher-genus algebraic curve. This case forced us to deal with some new techniques, such as separation of variables and finite-gap integration. One striking aspect was the role-reversal interplay between the discrete and continuous evolution of the separation variables. Another point of interest is the diagram of commuting canonical transformations between the original, separating and action-angle variables, which follows as a consequence of the transformation to action-angle variables on the one hand and the dynamical mapping on the other hand.

We finish by making a few general remarks on further ramifications. We remark that for the MHs which are constructed by the well-known BCH series, the precise relation between the integrability of mappings and the structure of the BCH series still remains mysterious. One may conjecture that a deeper understanding of the BCH series may arise from the convergence analysis of the expansion for the MH. We also remark that the interplay between discrete-time integrable systems and geometric integration occurs in a wider range of settings than that considered here. In both geometric integration in numerical analysis and quantum field theory, mathematical structures have been discovered and studied which provide insights into the BCH series, for instance rooted-trees expansions, Butcher groups and Hopf algebras. These insights have never been directly applied, to our knowledge, to the context of discrete integrable system. We may conjecture that techniques arising from the theory of discrete integrable system combined with the techniques arising from the areas of numerical analysis and Hopf algebras, may shed light on the structure of the BCH series.


\section*{Acknowledgement}

S.A.M. Alsallami would like to thank Umm Al-Qura University in Saudi Arabia and UKSACB in London for their financial support. F.W. Nijhoff is partially supported by EPSRC grant EP/I038683/1.

\appendix


\section{Lax Pairs}
\label{appendix:a}

Here we will go over the steps how to construct the Lax pair for the mappings of KdV type starting from the Lax representation of the lattice KdV equation (\ref{KdV eq}).
Recall the Lax representation of the lattice KdV equation,
$$\mathcal{L}={U}\,\mathrm{P}\,\widetilde{{U}}^{-1}, \qquad \mathcal{M}={U}\,\mathrm{Q}\,\widehat{{U}}^{-1}\ ,$$
in which
$${U}=\begin{pmatrix}1 & 0\\ u & 1\end{pmatrix},\quad \mathrm{P}=\begin{pmatrix}p & 1\\ k^2 & p\end{pmatrix},\quad \mathrm{Q}=\begin{pmatrix}q & 1\\ k^2 & q\end{pmatrix}.$$
As the initial value configuration that we have considered, the Lax matrices $\mathcal{L}$, $\mathcal{M}$ can be written as
$$
\mathcal{L}=\mathrm{U}_{2j}\,\mathrm{P}\,\mathrm{U}_{2j+1}^{-1}, \qquad \mathcal{M}=\mathrm{U}_{2j+1}\,\mathrm{Q}\,\mathrm{U}_{2j+2}^{-1}\ .
$$
	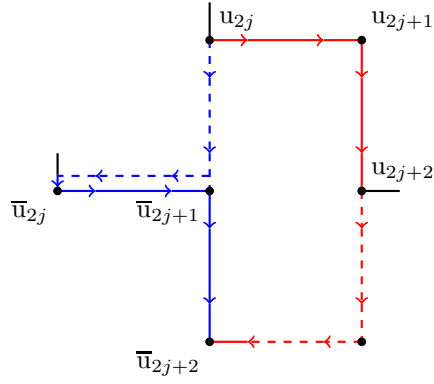
\begin{figure}[ht]
\centering
\begin{tikzpicture}[every node/.style={minimum size=0.1cm}]

\draw[thick] (2,0) -- (2.5,0);
\draw[thick] (0,2) -- (0,2.5);
\draw[thick] (-2,0.5) -- (-2,0.215);

\draw[color=red, thick,-] (1.5,2) -- (2,2);
\draw[color=red, thick,->] (0.5,2) -- (1.5,2);
\draw[color=red, thick,<-] (0.5,2) -- (0,2);
\draw[color=red, thick,-] (2,0) -- (2,0.5);
\draw[color=red, thick,<-] (2,0.5) -- (2,1.5);
\draw[color=red, thick,<-] (2,1.5) -- (2,2);
\draw[color=red, thick,-] (0.5,-2) -- (0,-2);
\draw[color=red, thick,dashed,->] (2,0) -- (2,-0.5);
\draw[color=red, thick,dashed,->] (2,-0.5) -- (2,-1.5);
\draw[color=red, thick,dashed,-] (2,-1.5) -- (2,-2);
\draw[color=red, thick,dashed,->] (1.5,-2) -- (0.5,-2);
\draw[color=red, thick,dashed,->] (2,-2) -- (1.5,-2);

\draw[color=blue, thick,<-] (0,-0.5) -- (0,0);
\draw[color=blue, thick,<-] (0,-1.5) -- (0,-0.5);
\draw[color=blue, thick,-] (0,-2) -- (0,-1.5);
\draw[color=blue, thick,<-] (-1.5,0) -- (-2,0);
\draw[color=blue, thick,<-] (-0.5,0) -- (-1.5,0);
\draw[color=blue, thick,-] (0,0) -- (-0.5,0);
\draw[color=blue, thick,<-] (-2,0.05) -- (-2,0.215);
\draw[color=blue, thick,-] (-2,0) -- (-2,0.05);
\draw[color=blue, thick,dashed,->] (0,2) -- (0,1.5);
\draw[color=blue, thick,dashed,->] (0,1.5) -- (0,0.5);
\draw[color=blue, thick,dashed,-] (0,0.5) -- (0,0);
\draw[color=blue, thick,dashed,->] (0,0.2) -- (-0.5,0.2);
\draw[color=blue, thick,dashed,->] (-0.5,0.2) -- (-1.5,0.2);
\draw[color=blue, thick,dashed,-] (-1.5,0.2) -- (-2,0.2);

\filldraw [fill=black] circle (0.05) 
          node[anchor= north east] {$\overline{\mathrm{u}}_{2j+1}$} ;
\filldraw [fill=black] (0,2) circle (0.05)
          node[anchor = south west] {$\mathrm{u}_{2j}$} ;
\filldraw [fill=black] (2,2) circle (0.05)
          node[anchor = south west] {$\mathrm{u}_{2j+1}$} ;         
\filldraw [fill=black] (2,0) circle (0.05)
          node[anchor = south west] {$\mathrm{u}_{2j+2}$} ;
\filldraw [fill=black] (2,-2) circle (0.05)
          node[anchor = south west] {} ;  
\filldraw [fill=black] (0,-2) circle (0.05)
          node[anchor = north east] {$\overline{\mathrm{u}}_{2j+2}$} ;          
\filldraw [fill=black] (-2,0) circle (0.05)
          node[anchor = north east] {$\overline{\mathrm{u}}_{2j}$} ;          
                
\end{tikzpicture}
	\caption{The Lax matrices for the mappings are derived by a consideration of the
compatibility of the two paths shown}
\label{figure of Lax matrices for the mappings}
\end{figure}
By consideration of the two paths as shown in figure~\ref{figure of Lax matrices for the mappings}, we have
\begin{align*}
	&\ol{\mathrm{U}}_{2j+3}\,\mathrm{P}^{-1}\,\ol{\mathrm{U}}_{2j+2}^{-1}\,
	{\mathrm{U}}_{2j+2}\,\mathrm{Q}\,\ol{\mathrm{U}}_{2j+3}^{-1}\,
	{\mathrm{U}}_{2j+1}\,\mathrm{Q}\,{\mathrm{U}}_{2j+2}^{-1}\,
	{\mathrm{U}}_{2j}\,\mathrm{P}\,{\mathrm{U}}_{2j+1}^{-1}\\
	&=\ol{\mathrm{U}}_{2j+1}\,\mathrm{Q}\,\ol{\mathrm{U}}_{2j+2}^{-1}\,
	\ol{\mathrm{U}}_{2j}\,\mathrm{P}\,\ol{\mathrm{U}}_{2j+1}^{-1}\,
	\ol{\mathrm{U}}_{2j+1}\,\mathrm{P}^{-1}\,\ol{\mathrm{U}}_{2j}^{-1}\,
	{\mathrm{U}}_{2j}\,\mathrm{Q}\,\ol{\mathrm{U}}_{2j+1}^{-1}\ ,\\
	\implies & \left(\mathrm{Q}\,\ol{\mathrm{U}}_{2j+2}^{-1}\,
	\ol{\mathrm{U}}_{2j}\,\mathrm{P}\,\ol{\mathrm{U}}_{2j+1}^{-1}\,\ol{\mathrm{U}}_{2j-1}\right)\left(\ol{\mathrm{U}}_{2j-1}^{-1}\,
	\ol{\mathrm{U}}_{2j+1}\,\mathrm{P}^{-1}\,\ol{\mathrm{U}}_{2j}^{-1}\,
	{\mathrm{U}}_{2j}\,\mathrm{Q}\,\ol{\mathrm{U}}_{2j+1}^{-1}\,{\mathrm{U}}_{2j-1}\right)\\
	&=\left(\ol{\mathrm{U}}_{2j+1}^{-1}\,\ol{\mathrm{U}}_{2j+3}\,\mathrm{P}^{-1}\,\ol{\mathrm{U}}_{2j+2}^{-1}\,{\mathrm{U}}_{2j+2}\,\mathrm{Q}\,\ol{\mathrm{U}}_{2j+3}^{-1}\,{\mathrm{U}}_{2j+1}\right)
	\left(\mathrm{Q}\,{\mathrm{U}}_{2j+2}^{-1}\,
	{\mathrm{U}}_{2j}\,\mathrm{P}\,{\mathrm{U}}_{2j+1}^{-1}\,{\mathrm{U}}_{2j-1}\right).
\end{align*}
Renaming the matrix products in the brackets, it is easy to show that we can write the latter equation in the form $$\ol{\mathbf{L}}_j\,\mathbf{M}_j=\mathbf{M}_{j+1}\,\mathbf{L}_j$$
in which
$$\mathbf{L}_j=\mathrm{Q}\,{\mathrm{U}}_{2j+2}^{-1}\,
	{\mathrm{U}}_{2j}\,\mathrm{P}\,{\mathrm{U}}_{2j+1}^{-1}\,{\mathrm{U}}_{2j-1},\quad \mathbf{M}_j=\ol{\mathrm{U}}_{2j-1}^{-1}\,
	\ol{\mathrm{U}}_{2j+1}\,\mathrm{P}^{-1}\,\ol{\mathrm{U}}_{2j}^{-1}\,
	\mathrm{U}_{2j}\,\mathrm{Q}\,\ol{\mathrm{U}}_{2j+1}^{-1}\,{\mathrm{U}}_{2j-1}\ ,$$
which can be used to construct the integrals of mappings. The Lax matrices $\mathbf{L}_j$ and $\mathbf{M}_j$ factorise as follows
$$\mathbf{L}_j=\mathsf{Q}\,L_j\,\mathsf{Q}^{-1}, \quad \mathbf{M}_j=\mathsf{Q}\,M_j\,\mathsf{Q}^{-1}, \qquad \mathsf{Q}=\begin{pmatrix}
		1 & 0\\
		q & 1
	\end{pmatrix},$$
in which the mappings of KdV type arises as the compatibility condition
$$\ol{L}_j\,M_j=M_{j+1}\,L_j\ ,$$
with the Lax matrices
$${L}_j=\begin{pmatrix}
		y_{j} & 1\\
		\ld & 0
	\end{pmatrix}\begin{pmatrix}
		x_{j} & 1\\
		\ld-\omega & 0
	\end{pmatrix}, \quad
{M}_j=\begin{pmatrix}
		{-\omega}/\ol{x}_j & 1\\
		\ld & -\ol{x}_j
	\end{pmatrix}\begin{pmatrix}
		x_{j}-{\omega}/{y_{j}} & 1\\
		\ld & 0
	\end{pmatrix}.$$
Using the same way as in the case of the lattice KdV, the Lax pair for the mappings of MKdV type can be constructed by using the relevant Lax representation given in factorized form as follows
$$\mathfrak{L}={V}^{-1}\,\mathrm{P}\,\widetilde{{V}}, \qquad \mathfrak{M}={V}^{-1}\,\mathrm{Q}\,\widehat{{V}}\ ,$$
in which
$${V}=\begin{pmatrix}1 & 0\\ 0 & v\end{pmatrix},\quad \mathrm{P}=\begin{pmatrix}p & 1\\ k^2 & p\end{pmatrix},\quad \mathrm{Q}=\begin{pmatrix}q & 1\\ k^2 & q\end{pmatrix}.$$


\section{Mapping Action}
\label{appendix:b}

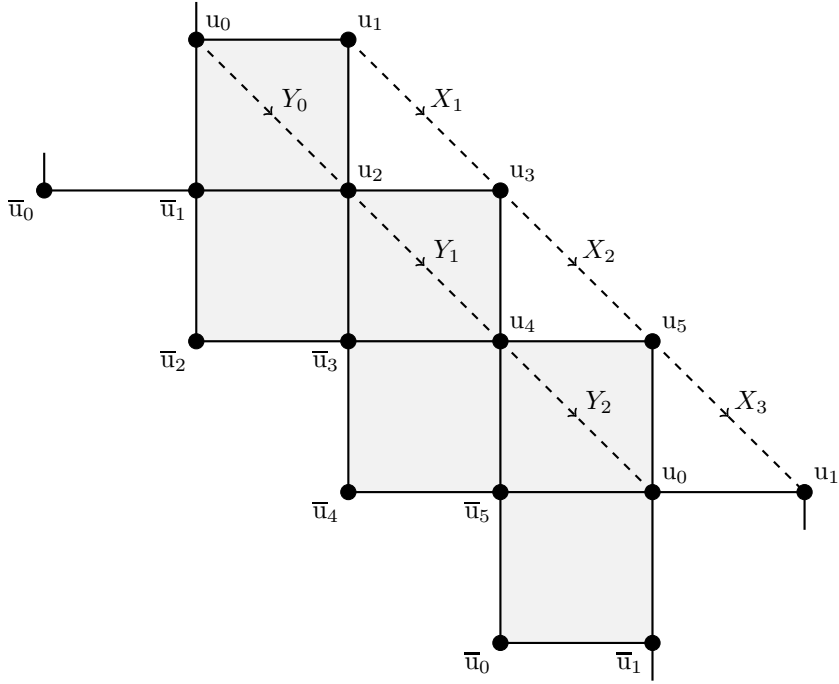
\begin{figure}[ht]
\centering
\begin{tikzpicture}[every node/.style={minimum size=0.2cm}]
\tikzstyle{solid node}=[circle,draw,inner sep=1,fill=black]

\draw[fill=black!5!white] 
    plot[smooth,samples=100,domain=0:2] (\x,{0}) -- 
    plot[smooth,samples=100,domain=2:0] (\x,{2});
\draw[fill=black!5!white] 
    plot[smooth,samples=100,domain=0:2] (\x,{0}) -- 
    plot[smooth,samples=100,domain=2:0] (\x,{-2});
\draw[fill=black!5!white] 
    plot[smooth,samples=100,domain=2:4] (\x,{0}) -- 
    plot[smooth,samples=100,domain=4:2] (\x,{-2});
\draw[fill=black!5!white] 
    plot[smooth,samples=100,domain=2:4] (\x,{-2}) -- 
    plot[smooth,samples=100,domain=4:2] (\x,{-4});
\draw[fill=black!5!white] 
    plot[smooth,samples=100,domain=4:6] (\x,{-2}) -- 
    plot[smooth,samples=100,domain=6:4] (\x,{-4});
\draw[fill=black!5!white] 
    plot[smooth,samples=100,domain=4:6] (\x,{-4}) -- 
    plot[smooth,samples=100,domain=6:4] (\x,{-6});
\draw[thick,-] (0,2) -- (0,2.5);
\draw[thick,-] (0,0) -- (2,0) node[above right=-0.3] {$\mathrm{u}_2$};
\draw[thick,-] (2,2) -- (0,2) node[above right=-0.3] {$\mathrm{u}_0$};
\draw[thick,-] (2,0) -- (2,2) node[above right=-0.3] {$\mathrm{u}_1$};
\draw[thick,-] (2,0) -- (4,0) node[above right=-0.3] {$\mathrm{u}_3$};
\draw[thick,-] (4,0) -- (4,-2) node[above right=-0.3] {$\mathrm{u}_4$};
\draw[thick,-] (4,-2) -- (6,-2);
\draw[thick,-] (6,-4) -- (6,-2) node[above right=-0.3] {$\mathrm{u}_5$};
\draw[thick,-] (0,2) -- (0,0) node[below left=-0.3] {$\overline{\mathrm{u}}_1$};
\draw[thick,-] (2,0) -- (2,-2) node[below left=-0.3] {$\overline{\mathrm{u}}_3$};
\draw[thick,-] (4,-2) -- (4,-4) node[below left=-0.3] {$\overline{\mathrm{u}}_5$};
\draw[thick,-] (6,-4) -- (8,-4);
\draw[thick,-] (6,-4) -- (6,-6) node[below left=-0.3] {$\ol{\mathrm{u}}_1$};
\draw[thick,-] (8,-4.5) -- (8,-4) node[above right=-0.3] {$\mathrm{u}_1$};
\draw[thick,-] (0,-2) -- (0,0);
\draw[thick,-] (2,-2) -- (0,-2) node[below left=-0.3] {$\overline{\mathrm{u}}_2$};
\draw[thick,-] (0,0) -- (-2,0) node[below left=-0.3] {$\overline{\mathrm{u}}_0$};
\draw[thick,-] (-2,0) -- (-2,0.5);
\draw[thick,-] (2,-2) -- (4,-2);
\draw[thick,-] (2,-4) -- (2,-2);
\draw[thick,-] (4,-4) -- (2,-4) node[below left=-0.3] {$\overline{\mathrm{u}}_4$};
\draw[thick,-] (4,-4) -- (4,-6);
\draw[thick,-] (4,-4) -- (6,-4) node[above right=-0.3] {$\mathrm{u}_0$};
\draw[thick,-] (6,-6) -- (4,-6) node[below left=-0.3] {$\overline{\mathrm{u}}_0$};
\draw[thick,-] (6,-6) -- (6,-6.5);
\draw[thick,dashed,->] (2,2) -- (3,1);
\draw[thick,dashed,-] (4,0) -- (3,1);
\draw[thick,dashed,->] (4,0) -- (5,-1);
\draw[thick,dashed,-] (6,-2) -- (5,-1);
\draw[thick,dashed,->] (0,2) -- (1,1);
\draw[thick,dashed,-] (2,0) -- (1,1);
\draw[thick,dashed,->] (2,0) -- (3,-1);
\draw[thick,dashed,-] (4,-2) -- (3,-1);
\draw[thick,dashed,->] (4,-2) -- (5,-3);
\draw[thick,dashed,-] (6,-4) -- (5,-3);
\draw[thick,dashed,->] (6,-2) -- (7,-3);
\draw[thick,dashed,-] (7,-3) -- (8,-4);
\node at (3.3,1.2) {$X_1$};
\node at (5.3,-0.8) {$X_2$};
\node at (7.3,-2.8) {$X_3$};
\node at (1.3,1.2) {$Y_0$};
\node at (3.3,-0.8) {$Y_1$};
\node at (5.3,-2.8) {$Y_2$};
\node[solid node] at (-2,0) {};
\node[solid node] at (0,0) {};
\node[solid node] at (2,0) {};
\node[solid node] at (4,0) {};
\node[solid node] at (0,2) {};
\node[solid node] at (2,2) {};
\node[solid node] at (0,-2) {};
\node[solid node] at (2,-2) {};
\node[solid node] at (4,-2) {};
\node[solid node] at (6,-2) {};
\node[solid node] at (2,-4) {};
\node[solid node] at (4,-4) {};
\node[solid node] at (6,-4) {};
\node[solid node] at (4,-6) {};
\node[solid node] at (6,-6) {};
\node[solid node] at (8,-4) {};
   \end{tikzpicture}
	\caption{Action of mapping over three periods}
\label{fig. of action for KdV g=2}
\end{figure}

We present here an action for the KdV lattice equation~(\ref{KdV eq}) for the case of genus two, as depicted in figure~\ref{fig. of action for KdV g=2}, namely
\begin{align*}
\mathcal{S}_{per}=\sum_t & \left\lbrack\mathrm{u}_{0}\left(\mathrm{u}_{1}-\ol{\mathrm{u}}_{1}\right)+\epsilon\,\delta\,\textrm{log}\left(\epsilon+\mathrm{u}_{0}-\mathrm{u}_{2}\right)+\ol{\mathrm{u}}_{1}\left(\mathrm{u}_{2}-\ol{\mathrm{u}}_{2}\right)+\epsilon\,\delta\,\textrm{log}\left(\epsilon+\ol{\mathrm{u}}_{1}-\ol{\mathrm{u}}_{3}\right)\right. \\
&+\mathrm{u}_{2}\left(\mathrm{u}_{3}-\ol{\mathrm{u}}_{3}\right)+\epsilon\,\delta\,\textrm{log}\left(\epsilon+\mathrm{u}_{2}-\mathrm{u}_{4}\right)+\ol{\mathrm{u}}_{3}\left(\mathrm{u}_{4}-\ol{\mathrm{u}}_{4}\right)+\epsilon\,\delta\,\textrm{log}\left(\epsilon+\ol{\mathrm{u}}_{3}-\ol{\mathrm{u}}_{5}\right) \\
&\left. +\,\mathrm{u}_{4}\left(\mathrm{u}_{5}-\ol{\mathrm{u}}_{5}\right)+\epsilon\,\delta\,\textrm{log}\left(\epsilon+\mathrm{u}_{4}-\mathrm{u}_{0}\right)+\ol{\mathrm{u}}_{5}\left(\mathrm{u}_{0}-\ol{\mathrm{u}}_{0}\right)+\epsilon\,\delta\,\textrm{log}\left(\epsilon+\ol{\mathrm{u}}_{5}-\ol{\mathrm{u}}_{1}\right)\right\rbrack ,
\end{align*}
in which $\mathrm{u}=\mathrm{u}(t), \ol{\mathrm{u}}=\mathrm{u}(t+1)$ and $\mathcal{S}_{per}$ is $\mathcal{S}$~(\ref{eq:action for KdV}) modulo periodicity.
Using the reduced variables
	$$X_j\equiv \mathrm{u}_{{2j}+1}-\mathrm{u}_{{2j}-1},\qquad Y_j\equiv \mathrm{u}_{{2j}+2}-\mathrm{u}_{2j}\ ,$$
the latter equation can be reduced to the following form
\be
\begin{split}
	\mathcal{S}_{per}=\sum_t &\left\lbrack Y_1\,\ol{X}_1-Y_1\,X_1-Y_2\,\ol{X}_3+Y_2\,X_3\right. \\
	&+\epsilon\,\delta\,\textrm{log}\,(\epsilon-Y_1)(\epsilon-Y_2)(\epsilon+Y_1+Y_2) \\
	&\left. +\epsilon\,\delta\,\textrm{log}\,(\epsilon-\ol{X}_1)(\epsilon-\ol{X}_3)(\epsilon+\ol{X}_1+\ol{X}_3)\right\rbrack ,
\end{split}
\label{action eq KdV g=2}
\ee
where $X_3=-(X_1+X_2)$ and $Y_0=-(Y_1+Y_2)$ using the periodicity constraints.
The Euler-Lagrange equations for (\ref{action eq KdV g=2}), which are obtained by variation of $\mathcal{S}$ with respect to the variables $X_1, X_3, Y_1, Y_2$, i.e.\
\begin{align*}
	\frac{\updelta\mathcal{S}}{\updelta X_1}&=\underline{Y_1}-Y_1-\frac{\epsilon\,\delta}{\epsilon-X_1}+\frac{\epsilon\,\delta}{\epsilon+X_1+X_3}=0\ ,\\
	\frac{\updelta\mathcal{S}}{\updelta X_3}&=Y_2-\underline{Y_2}-\frac{\epsilon\,\delta}{\epsilon-X_3}+\frac{\epsilon\,\delta}{\epsilon+X_1+X_3}=0\ ,\\
	\frac{\updelta\mathcal{S}}{\updelta Y_1}&=\ol{X}_1-X_1-\frac{\epsilon\,\delta}{\epsilon-Y_1}+\frac{\epsilon\,\delta}{\epsilon+Y_1+Y_2}=0\ ,\\
	\frac{\updelta\mathcal{S}}{\updelta Y_2}&=X_3-\ol{X}_3-\frac{\epsilon\,\delta}{\epsilon-Y_2}+\frac{\epsilon\,\delta}{\epsilon+Y_1+Y_2}=0\ ,
\end{align*}
lead to the equations
\bse
\label{eqs for variation KdV g=2}
\begin{align}
\ol{Y}_1-Y_1&=\frac{\epsilon\,\delta}{\epsilon+\ol{X}_1+\ol{X}_3}-\frac{\epsilon\,\delta}{\epsilon-\ol{X}_1},\quad \ol{Y}_2-Y_2=\frac{\epsilon\,\delta}{\epsilon-\ol{X}_3}-\frac{\epsilon\,\delta}{\epsilon+\ol{X}_1+\ol{X}_3}\ , \\
\ol{X}_1-X_1&=\frac{\epsilon\,\delta}{\epsilon-{Y}_1}-\frac{\epsilon\,\delta}{\epsilon+Y_1+Y_2}\, \  \ ,\quad \ol{X}_3-X_3=\frac{\epsilon\,\delta}{\epsilon+Y_1+Y_2}-\frac{\epsilon\,\delta}{\epsilon-{Y}_2} \quad ,
\end{align}
\ese
which is the integrable mapping~(\ref{mapping of KdV (P=3)in X and Y}). Now, our aim is to perform a Legendre transformation to establish a Hamiltonian structure for the mapping. A discussion of performing such a discrete-time Legendre transformation can be found in~\cite{bruschi1991integrable}.

Consider a Lagrangian $\mathscr{L}=\mathscr{L}(X_1, \ol{X}_1, X_3, \ol{X}_3, Y_1, Y_2)$, such that the Euler-Lagrange equations
$$\ol{\frac{{\partial\mathscr{L}}}{\partial X_1}}+\frac{\partial\mathscr{L}}{\partial \ol{X}_1}=0, \qquad \ol{\frac{{\partial\mathscr{L}}}{\partial X_3}}+\frac{\partial\mathscr{L}}{\partial \ol{X}_3}=0\ ,$$
correspond to the mapping under consideration. The discrete-time Hamiltonian is obtained through the following Legendre transformation
\begin{align*}
{H}(Y_1, Y_2, \ol{X}_1, \ol{X}_3)=&\,Y_1\,X_1-Y_1\,\ol{X}_1+Y_2\,\ol{X}_3-Y_2\,X_3+\mathscr{L}\ .
\end{align*}
Now, we have
\begin{align*}
	&\frac{\partial {H}}{\partial \ol{X}_1}\,\updelta \ol{X}_1+\frac{\partial {H}}{\partial \ol{X}_3}\,\updelta \ol{X}_3+\frac{\partial {H}}{\partial {Y}_1}\,\updelta {Y}_1+\frac{\partial {H}}{\partial {Y}_2}\,\updelta {Y}_2\\
	&=Y_1\,\updelta {X}_1+X_1\,\updelta {Y}_1-Y_1\,\updelta \ol{X}_1-\ol{X}_1\,\updelta {Y}_1
	+\ol{X}_3\,\updelta {Y}_2+Y_2\,\updelta \ol{X}_3-Y_2\,\updelta {X}_3-X_3\,\updelta {Y}_2 \\
	&~\  \ +\frac{\partial \mathscr{L}}{\partial {X}_1}\,\updelta {X}_1+\frac{\partial \mathscr{L}}{\partial \ol{X}_1}\,\updelta \ol{X}_1+\frac{\partial \mathscr{L}}{\partial {X}_3}\,\updelta {X}_3+\frac{\partial \mathscr{L}}{\partial \ol{X}_3}\,\updelta \ol{X}_3+\frac{\partial \mathscr{L}}{\partial {Y}_1}\,\updelta {Y}_1+\frac{\partial \mathscr{L}}{\partial {Y}_2}\,\updelta {Y}_2\ .
\end{align*}
Thus, one obtains a set of equations
\begin{align*}
	\frac{\partial {H}}{\partial \ol{X}_1}&=-Y_1+\frac{\partial \mathscr{L}}{\partial \ol{X}_1}\quad~\  \ , \qquad \frac{\partial {H}}{\partial \ol{X}_3}=Y_2+\frac{\partial \mathscr{L}}{\partial \ol{X}_3}~\qquad , \\
	\frac{\partial {H}}{\partial {Y}_1}&=X_1-\ol{X}_1+\frac{\partial \mathscr{L}}{\partial {Y}_1}\ , \qquad\, \frac{\partial {H}}{\partial {Y}_2}=\ol{X}_3-X_3+\frac{\partial \mathscr{L}}{\partial {Y}_2}\ , \\
	0&={Y}_1+\frac{\partial \mathscr{L}}{\partial {X}_1}\qquad~ ,\  \ \quad\qquad 0=-{Y}_2+\frac{\partial \mathscr{L}}{\partial {X}_3}\quad~\ \ ,
\end{align*}
and consequently we get
\begin{align*}
	\frac{\partial {H}}{\partial \ol{X}_1}&=\ol{Y}_1-Y_1\, , \qquad \frac{\partial {H}}{\partial \ol{X}_3}=Y_2-\ol{Y}_2\ \, , \\
	\frac{\partial {H}}{\partial {Y}_1}&=X_1-\ol{X}_1, \qquad \frac{\partial {H}}{\partial {Y}_2}=\ol{X}_3-X_3\ .
\end{align*}
By identifying
$$X_1:= p_1,\qquad X_3:= -p_2,\qquad Y_1:= q_1,\qquad Y_2:= q_2\ ,$$
one obtains the discrete-time Hamiltonian in terms of variables $q_1, q_2, \ol{p}_1, \ol{p}_2$
\begin{equation*}
\begin{split}
    {H}(q_1,q_2,\overline{p}_1,\overline{p}_2)=&\,\epsilon \,\delta\, \textrm{log}\,(\epsilon+\overline{p}_1-\overline{p}_2)(\epsilon-\overline{p}_1)(\epsilon+\overline{p}_2)\\
    &+\epsilon \,\delta\, \textrm{log}\,(\epsilon+q_1+q_2)(\epsilon-q_1)(\epsilon-q_2)\ ,
    \end{split}
\end{equation*}
which acts as the generating functional for the mapping, i.e. one has the discrete-time Hamilton equations
\begin{equation*}
\begin{split}
	&\overline{p}_1-p_1=-\frac{\partial {H}}{\partial q_1}, \qquad   \ \overline{q}_1-q_1=\frac{\partial {H}}{\partial \overline{p}_1}\ , \\
	&\overline{p}_2-p_2=-\frac{\partial {H}}{\partial q_2}, \qquad   \ \overline{q}_2-q_2=\frac{\partial {H}}{\partial \overline{p}_2}\ ,
	\end{split}
\end{equation*}
leading to the equations of the mapping under consideration. In terms of the conjugate variables $q_1,q_2,p_1,p_2$ we have the standard symplectic structure
$$
\Omega=dq_1 \wedge dp_1+dq_2 \wedge dp_2\ ,
$$
leading to standard Poisson brackets
$$\{q_i,q_{j}\}=\{p_i,p_{j}\}=0,\quad \{q_i,p_{j}\}= \delta_{ij},\qquad i,j=1,2 ,$$ and the mapping is in fact a canonical transformation, i.e.\ $\ol{\Omega}=\Omega$.


\section{Proof of Equation (\ref{Poisson structure for monodromy matrix})}
\label{appendix:c}

In this appendix we prove the relation (\ref{Poisson structure for monodromy matrix}). Let us first recall the discrete version of the non-ultralocal Poisson bracket structure from~\cite{nijhoff1992integrable}, i.e.
\begin{equation*}
\{\overset{1}{L}_n,\overset{2}{L}_m\}=-\dd_{n,m+1}\,\overset{1}{L}_n\, s_{12}^{+}\, \overset{2}{L}_m+\dd_{n+1,m}\, \overset{2}{L}_m\, s_{12}^{-}\, \overset{1}{L}_n
+\dd_{n,m}\,(r_{12}^{+}\,\overset{1}{L}_n\,\overset{2}{L}_m-\overset{1}{L}_n\,\overset{2}{L}_m\,r_{12}^{-})\ .
\end{equation*}
Using this equation together with equation~(\ref{eq:monodrmy matrix by Lax (SoV sec)}) we can establish
\begin{equation*}
	\begin{split}
		\{\overset{1}{T},\overset{2}{L}_n\}\,&{=\{\overset{1}{L}_P\,\overset{1}{L}_{P-1}\cdots\overset{1}{L}_1,\,\overset{2}{L}_n\}}\\
		&{\begin{split}
			 =&\,\{\overset{1}{L}_{P},\,\overset{2}{L}_{n}\}\,\overset{1}{L}_{P-1}\cdots\overset{1}{L}_{1}+\overset{1}{L}_{P}\cdots\overset{1}{L}_{2}\, \{\overset{1}{L}_{1},\,\overset{2}{L}_{n}\}\\
			 &+\sum_{j=2}^{P-1} \overset{1}{L}_{P}\cdots\overset{1}{L}_{j+1}\,\{\overset{1}{L}_{j},\,\overset{2}{L}_{n}\}\, \overset{1}{L}_{j-1}\cdots\overset{1}{L}_{1} \end{split}}\\
		&{\begin{split}=&\,\{\overset{1}{L}_{P},\,\overset{2}{L}_{n}\}\,\overset{1}{L}_{P-1}\cdots\overset{1}{L}_{1}
		+\overset{1}{L}_{P}\cdots\overset{1}{L}_{n+2}\left(-\overset{1}{L}_{n+1}\,s_{12}^{+}\, \overset{2}{L}_{n}\right)\overset{1}{L}_{n}\cdots\overset{1}{L}_{1}\\
		&+\overset{1}{L}_{P}\cdots\overset{1}{L}_{n+1}\left(r_{12}^{+}\,\overset{1}{L}_{n}\,\overset{2}{L}_{n}-\overset{2}{L}_{n}\,\overset{1}{L}_{n}\,r_{12}^{-}\right)\overset{1}{L}_{n-1}\cdots\overset{1}{L}_{1}\\
		&+\overset{1}{L}_{P}\cdots\overset{1}{L}_{n}\left(\overset{2}{L}_{n}\,s_{12}^{-}\, \overset{1}{L}_{n-1}\right)\overset{1}{L}_{n-2}\cdots\overset{1}{L}_{1}+\overset{1}{L}_{P}\cdots\overset{1}{L}_{2}\, \{\overset{1}{L}_{1},\,\overset{2}{L}_{n}\}\ .\end{split}}
	\end{split}
	\end{equation*}
At this point it is useful to establish the commutation relations between the monodromy matrices $T$ and $T_n^{+}, T_n^{-}$, using the fundamental commutation relations of the matrices $L_n$.
Introduce the following decomposition of the monodromy matrix (\ref{eq:monodrmy matrix by Lax (SoV sec)})$$T=T_n^{+}\cdot T_{n-1}^{-}\ ,$$in which
	$$T_n^{+}(\ld)=\prod_{j=n}^{P}L_j(\ld),\quad T_{n-1}^{-}(\ld)=\prod_{j=1}^{n-1}L_j(\ld)\ .$$
Thus, we obtain
\begin{equation*}
\begin{split}
	\{\overset{1}{T},\overset{2}{L}_n\}=&\,\{\overset{1}{L}_P,\overset{2}{L}_n\}\,\overset{1}{T}{_{P-1}^{-}}+\overset{2}{T}{_{2}^{+}}\,\{\overset{1}{L}_{1},\,\overset{2}{L}_{n}\}
	+\overset{2}{L}_{n}\,\overset{1}{T}{_{n}^{+}}\left(s_{12}^{-}-r_{12}^{-}\right)\overset{1}{T}{_{n-1}^{-}}\\
	&-\overset{1}{T}{_{n+1}^{+}}\left(s_{12}^{+}-r_{12}^{+}\right)\overset{1}{T}{_{n}^{-}}\,\overset{2}{L}_{n}	\end{split}
\end{equation*}
Finally, we have that
	\begin{equation*}
		\begin{split}
		\{\overset{1}{T},\overset{2}{T}\}\,&{=\sum_{n=1}^{P} \overset{2}{L}_{P}\cdots\overset{2}{L}_{n+1}\,\{\overset{1}{T},\,\overset{2}{L}_{n}\}\, \overset{2}{L}_{n-1}\cdots\overset{2}{L}_{1}}\\
		&{=\{\overset{1}{T},\overset{2}{L}_P\}\,\overset{2}{T}{_{P-1}^{-}}+\overset{2}{T}{_{2}^{+}}\,\{\overset{1}{T},\overset{2}{L}_1\}+\sum_{n=2}^{P-1} \overset{2}{T}{_{n+1}^{+}}\,\{\overset{1}{T},\,\overset{2}{L}_{n}\}\, \overset{2}{T}{_{n-1}^{-}}}\\
		&{\begin{split}=&\,\{\overset{1}{T},\overset{2}{L}_P\}\,\overset{2}{T}{_{P-1}^{-}}+\sum_{n=2}^{P-1} \overset{2}{T}{_{n+1}^{+}}\,\{\overset{1}{L}_P,\,\overset{2}{L}_{n}\}\,\overset{1}{T}{_{P-1}^{-}}\, \overset{2}{T}{_{n-1}^{-}}\\
		&+\sum_{n=2}^{P-1}\left(\overset{2}{T}{_{n}^{+}}\,\overset{1}{T}{_{n}^{+}}\left(s_{12}^{-}-r_{12}^{-}\right)\overset{1}{T}{_{n-1}^{-}}\,\overset{2}{T}{_{n-1}^{-}}-\overset{2}{T}{_{n+1}^{+}}\, \overset{1}{T}{_{n+1}^{+}}\left(s_{12}^{+}-r_{12}^{+}\right)\overset{1}{T}{_{n}^{-}}\,\overset{2}{T}{_{n}^{-}}\right)\\
		&+\overset{2}{T}{_{2}^{+}}\,\{\overset{1}{T},\overset{2}{L}_1\}+\sum_{n=2}^{P-1} \overset{2}{T}{_{n+1}^{+}}\,\overset{1}{T}{_{2}^{+}}\{\overset{1}{L}_{1},\,\overset{2}{L}_{n}\}\,\overset{2}{T}{_{n-1}^{-}}\end{split}}\\
		&{\begin{split}=&\,\,r_{12}^{+}\,\overset{1}{T}\,\overset{2}{T}-\overset{2}{T}\,\overset{1}{T}\, r_{12}^{-}-\overset{1}{T}{_{2}^{+}}\,\overset{1}{T}{_{1}^{-}}\,s_{12}^{+} \, \overset{2}{T}{_{P}^{+}}\,\overset{2}{T}{_{P-1}^{-}}+\overset{2}{T}{_{P}^{+}}\,\overset{1}{T}{_{n}^{+}}\left(s_{12}^{-}-r_{12}^{-} \right) \overset{1}{T}{_{n-1}^{-}}\,\overset{2}{T}{_{P-1}^{-}}\\
		&-\overset{1}{T}{_{n+1}^{+}}\left(s_{12}^{+}-r_{12}^{+} \right) \overset{1}{T}{_{n}^{-}}\,\overset{2}{T}{_{P}^{-}}-\overset{2}{T}{_{P}^{+}}\,\overset{1}{T}{_{P}^{+}}\, s_{12}^{+} \, \overset{1}{T}{_{P-1}^{-}}\,\overset{2}{T}{_{P-1}^{-}}\\
		&+\overset{2}{T}{_{2}^{+}}\,\overset{1}{T}{_{2}^{+}}\left(s_{12}^{-}-r_{12}^{-} \right) \overset{1}{T}{_{1}^{-}}\,\overset{2}{T}{_{1}^{-}}-\overset{2}{T}{_{P}^{+}}\,\overset{1}{T}{_{P}^{+}}\left(s_{12}^{+}-r_{12}^{+} \right) \overset{1}{T}{_{P-1}^{-}}\,\overset{2}{T}{_{P-1}^{-}}\\
		&+\overset{2}{T}{_{2}^{+}}\,\overset{2}{T}{_{1}^{-}}\,s_{12}^{-} \, \overset{1}{T}{_{P}^{+}}\,\overset{1}{T}{_{P-1}^{-}}+r_{12}^{+}\,\overset{1}{T}\,\overset{2}{T}-\overset{2}{T}\,\overset{1}{T}\, r_{12}^{-}+\overset{2}{T}{_{1}^{+}}\,\overset{1}{T}{_{n}^{+}}\left(s_{12}^{-}-r_{12}^{-} \right) \overset{1}{T}{_{n-1}^{-}}\\
		&-\overset{2}{T}{_{2}^{+}}\,\overset{1}{T}{_{n+1}^{+}}\left(s_{12}^{+}-r_{12}^{+}\right)\overset{1}{T}{_{n}^{-}}\,\overset{2}{T}{_{1}^{-}}+\overset{2}{T}{_{2}^{+}}\, \overset{1}{T}{_{2}^{+}}\, s_{12}^{-}\, \overset{1}{T}{_{1}^{-}}\,\overset{2}{T}{_{1}^{-}}\end{split}}\\
		&{=r_{12}^{+}\,\overset{1}{T}\,\overset{2}{T}-\overset{2}{T}\,\overset{1}{T}\,r_{12}^{-}-\overset{1}{T}\, s_{12}^{+}\, \overset{2}{T}+\overset{2}{T}\, s_{12}^{-}\, \overset{1}{T}}\ .
		\end{split}
	\end{equation*}


\section{Discrete Dubrovin Equations}
\label{appendix:d}

In this appendix we will show how to derive the discrete Dubrovin equations for the mappings of KdV type starting from the the discrete-time evolution (\ref{map T (SoV)}).
The auxiliary spectrum $B(\ld)$ (\ref{eq:auxiliary spectrum B}) leads to the expressions for $B_0/B_g,\dots,B_{g-1}/B_g$ as elementary symmetric functions of the zeroes $\mu_1,\dots,\mu_g$. Similarly, the coefficients of $C(\ld)$ are symmetric functions of its zeroes $\gamma_1,\dots,\gamma_g$. From the fact that
\be \label{eq: disc R} 
(A-D)(\ld)=\kp\sqrt{R(\ld)-4\,B(\ld)C(\ld)}\   ,   
\ee
where the $\kappa$ denotes the sign $\kp=\pm$ corresponding to the choice of sheet of the Riemann surface, subject to the condition $\ol{\kp}=\kp$, and taking $\ld=\mu_i$ in (\ref{eq: disc R}) we obtain the system
\begin{equation*}\label{eq:system} 
\left( \begin{array}{c} 
\mathcal{I}_1 - 2D_0 \\ \vdots \\ \mathcal{I}_g - 2D_{g-1} \end{array}\right)  = 
\left(\begin{array}{ccc}
\mu_1 & \cdots & \mu_1^g \\ 
\vdots & & \vdots \\ 
\mu_g & \cdots &\mu_g^g  
\end{array} \right)^{-1}\cdot 
\left( \begin{array}{c} 
\kp_1\sqrt{R(\mu_1)}- \mathcal{I}_0 \\ \vdots \\ 
\kp_g\sqrt{R(\mu_g)}- \mathcal{I}_0 \end{array} \right) .
\end{equation*}
Thus, we obtain
\be\label{eq:1 proof of Dubv}
\mathcal{I}_j-2\,D_{j-1}=\left\lbrack{\mathscr M}^{-1}\left(\boldsymbol{\kp}\sqrt{R({\boldsymbol{\mu}})}-\mathcal{I}_0\,{\boldsymbol e}\right)\right\rbrack_j, \qquad j=1,\dots,g\ ,
\ee
where $\mathscr M$ is the VanderMonde matrix given in equation~(\ref{eq: VanderMonde matrix}).
From the Lax equation for the map (\ref{map T (SoV)}) we have the following discrete relations for its entries
\bse\label{eq:relations for T map entries}
\begin{align}
	&\left\lbrack\ld-\frac{\omega}{\ol{x}}\left(x-\frac{\omega}{y}\right)\right\rbrack A-\frac{\omega}{\ol{x}}\,C=\ol{A}\left\lbrack\ld-\frac{\omega}{\ol{x}}\left(x-\frac{\omega}{y}\right)\right\rbrack+\ld\,\ol{B}\left(x-\ol{x}-\frac{\omega}{y}\right),\label{eq: T map entry11} \\
	&\left\lbrack\ld-\frac{\omega}{\ol{x}}\left(x-\frac{\omega}{y}\right)\right\rbrack B-\frac{\omega}{\ol{x}}\,D=\ld\,\ol{B}-\frac{\omega}{\ol{x}}\,\ol{A}\ ,\label{eq: T map entry12} \\
	&\ld\left(x-\ol{x}-\frac{\omega}{y}\right)A+\ld\,C=\ol{C}\left\lbrack\ld-\frac{\omega}{\ol{x}}\left(x-\frac{\omega}{y}\right)\right\rbrack+\ld\,\ol{D}\left(x-\ol{x}-\frac{\omega}{y}\right),\label{eq: T map entry21} \\
	&\ld\left(x-\ol{x}-\frac{\omega}{y}\right)B+\ld\,D=\ld\,\ol{D}-\frac{\omega}{\ol{x}}\,\ol{C}\ .\label{eq: T map entry22}
\end{align}\ese
From (\ref{eq:relations for T map entries}) we establish that the top coefficients $B_g$ and $C_g$ can be taken to be equal and constant.
Expanding equation~(\ref{eq: T map entry11}) in powers of $\ld$ we are lead to the following coupled equations
\begin{align}
 \ol{A}_0 = A_0=\mathcal{I}_0\qquad ,\qquad
  A_g-\frac{\omega}{\ol{x}}\,C_g=\ol{A}_g+\ol{B}_g\left(x-\ol{x}-\frac{\omega}{y}\right) , \label{eq:1,2 from entry11}
\end{align}
while expanding equation~(\ref{eq: T map entry12}) we are lead to the following set of equations
\bse\label{eq:Exp from entry12}
\begin{align}
    &\ol{A}_0 = B_0\left(x-\frac{\omega}{y}\right)\qquad ,\qquad
  \ol{B}_g = B_g\ ,\label{eq:1,2 from entry12} \\
  &\frac{\omega}{\ol{x}}\left(\ol{A}_j-D_{j-1}\right) = \ol{B}_{j-1}-B_{j-1}+\frac{\omega}{\ol{x}}\left(x-\frac{\omega}{y}\right)B_j\ , \qquad j=1,\dots, g\ .\label{eq:3 from entry12}
\end{align}\ese
Using the equations~(\ref{eq:1,2 from entry11}) and (\ref{eq:1,2 from entry12}) one obtains
\be
B_g\,\ol{x}-\frac{\omega}{\ol{x}}\,C_g=D_{g-1}-\ol{D}_{g-1}+\frac{B_g}{B_0}\,\ol{A}_0\ .
\label{eq:D6 Dub eq}
\ee
By using the relation~(\ref{eq:1 proof of Dubv}) and the fact that
$$B(\ld) = B_g\prod_{j=1}^g (\ld-\mu_j)=\sum_{j=0}^g \ld^j\,B_j=B_g\sum_{j=0}^g\ld^j(-1)^{g-j}S_{g-j}(\boldsymbol{\mu})\ ,$$
together with the relations~(\ref{eq:3 from entry12}) and (\ref{eq:D6 Dub eq}) we obtain a coupled system of set of first-order difference equations for the $\mu_j$, respectively, namely
\bse\label{eq:coupled system, pf of Dub eq}
\begin{align} \label{eq1:coupled system, pf of Dub eq} \nonumber
&\left\lbrack{\mathscr M}^{-1}\left(\boldsymbol{\kp}\sqrt{R({\boldsymbol{\mu}})}-\mathcal{I}_0\,{\boldsymbol e}\right)\right\rbrack_j
 +\left\lbrack\xbar{\mathscr M}^{-1}\left(\ol{\boldsymbol{\kp}}\sqrt{R(\ol{\boldsymbol{\mu}})}-\mathcal{I}_0\,{\boldsymbol e}\right)\right\rbrack_j \\
&=\frac{2\,\mathcal{I}_0\,S_{g-j}(\boldsymbol{\mu})}{(-1)^{j}\prod_{i=1}^g\mu_i}+2\,B_g\,\frac{\ol{x}}{\omega}\,(-1)^{g-j+1}\left[ S_{g-j+1}(\ol{\boldsymbol{\mu}})-S_{g-j+1}(\boldsymbol{\mu})\right] , \\ \nonumber
&\left[{\mathscr M}^{-1}\left(\boldsymbol{\kp}\sqrt{R({\boldsymbol{\mu}})}-\mathcal{I}_0\,{\boldsymbol e}\right)\right]_g -\left[\xbar{\mathscr M}^{-1}\left(\ol{\boldsymbol{\kp}}\sqrt{R(\ol{\boldsymbol{\mu}})}-\mathcal{I}_0\,{\boldsymbol e}\right)\right]_g \\
&=\frac{2\,\mathcal{I}_0}{(-1)^{g}\prod_{i=1}^g\mu_i}-2\,B_g\,\ol{x}+2\,C_g\,\frac{\omega}{\ol{x}}\ ,
\label{eq2:coupled system, pf of Dub eq}
\end{align}\ese
where $j=1,\dots,g$.


\section{Baker-Akhiezer Eigenfunction}
\label{appendix:e}

Here we look at the asymptotic form of the Bloch eigenfunctions $\phi(\ld)$, i.e.~(\ref{eq:phi vector}), of the spectral problem. In fact, we have for these eigenfunctions the following properties
$$ T\phi^{\pm}=\upeta^{\pm}\phi^{\pm}, \qquad \ol{\phi}^{\pm}=M{\phi}^{\pm}\ ,$$
where $\upeta^{\pm}$ are the two different solutions of the hyperelliptic spectral curve, and the signs ${\pm}$ correspond to the solutions on the two sheets of $\Gamma$. We have the following expressions for the Bloch eigenfunctions
\be
\begin{pmatrix} 1 \\[0.1cm] \frac{\upeta^{\pm}-A}{B}\end{pmatrix}\varphi^{\pm}
	=\begin{pmatrix} 1 \\[0.1cm] \frac{C}{\upeta^{\pm}-D}\end{pmatrix}\varphi^{\pm}\ .
\label{eq:expres. Bloch eigenfun}
\ee
We also have
\be
\ol{\varphi}^{\pm}\begin{pmatrix} 1 \\[0.1cm] \frac{\ol{\upeta}^{\pm}-\ol{A}}{\ol{B}}\end{pmatrix}
={\varphi}^{\pm}\begin{pmatrix}
		\frac{\omega}{\ol{x}}\left(\frac{\omega}{y}-x\right)+\ld-\frac{\omega}{\ol{x}}\left(\frac{\upeta^{\pm}-A}{B}\right)\\
		\ld\left(x-\ol{x}-\frac{\omega}{y}\right)+\ld\left(\frac{\upeta^{\pm}-A}{B}\right)
\end{pmatrix}.
\label{eq: Baker eqE2}
\ee
The scalar function $\varphi^{\pm}$ solves
\be
\frac{\ol{\varphi}^{\pm}}{{\varphi}^{\pm}}=\frac{\omega}{\ol{x}}\left(\frac{\ol{D}-\upeta^{\pm}}{B}\right)+\ld\,\frac{\ol{B}}{B}\ ,
\label{eq: Baker eqE3}
\ee
thus we have
\be
\frac{\ol{\varphi}^{+}}{\varphi^{+}}\frac{\ol{\varphi}^{-}}{\varphi^{-}}=\ld\,(\ld-\omega)\,\frac{\ol{B}}{B}=\frac{(\ld-\omega)\,\ol{B}}{(1/\ld)\,B}\ .
\label{eq: Baker eqE4}
\ee
Note that~(\ref{eq: Baker eqE4}) has pole at $\infty$ and zero at $\ld=\omega$.
We point out that the $\mu_j$ arise as the poles, and the $\ol{\mu}_j$ as the zeroes of a transition factor ${\ol{\varphi}}/{{\varphi}}$ where $\varphi$ is the relevant Baker-Akhiezer function.

As a consequence of Abel's theorem, cf.~\cite{buchstaber1996hyperelliptic}, we have the following relation
\be
\sum_{j=1}^{g}\int_{(\mu_j,\eta_j)}^{(\ol{\mu}_j,\ol{\eta}_j)}\omega_k+\int_{\infty}^{(\omega,\upeta(\omega))}\omega_k=0\qquad\left(\textrm{mod}~\Lambda_g\right) ,
\label{eq: Baker eqE5}
\ee
which provides us with the solution of the discrete Dubrovin equations. In equation~(\ref{eq: Baker eqE5}) $\omega_k$ is the normalized differential of the first kind~(Abelian differential) and $\omega$ denotes the parameter value given in~(\ref{eq:discriminant(KdV)}), whereas $\Lambda_g$ is a period lattice of the associated Riemann surface.


\section*{References}


\end{document}